\documentclass[10pt]{article}
\usepackage{amssymb}
\newtheorem{theorem}{Theorem}

\newtheorem{corollary}[theorem]{Corollary}

\newtheorem{proposition}{Proposition}

\begin{document}
\title{Inverse  problem by Cauchy data on arbitrary
subboundary for system of elliptic equations}

\author{
O.~Yu.~Imanuvilov\thanks{ Department
of Mathematics, Colorado State University, 101 Weber Building, Fort
Collins, CO 80523-1874, USA, E-mail: oleg@math.colostate.edu}\, and
M.~Yamamoto\thanks{ Department of Mathematical Sciences, The University
of Tokyo, Komaba Meguro Tokyo 153 Japan,
E-mail:myama@ms.u-tokyo.ac.jp}
}

\date{}
\maketitle

\begin{abstract}
We consider an inverse problem of determining coefficient matrices
in an $N$-system of second-order elliptic equations in a bounded two
dimensional domain by a set of Cauchy data on arbitrary subboundary.
The main result of the article is as follows: If two systems of
elliptic operators generate the same set of partial Cauchy data on
an arbitrary subboundary, then the coefficient matrices of the
first-order and zero-order terms satisfy the prescribed system of
first-order partial differential equations.  The main result implies
the uniqueness of any two coefficient matrices provided that the one
remaining matrix among the three coefficient matrices is known.
\end{abstract}

\section{Introduction}\label{sec1}
Let $\Omega$ be a bounded domain in $\Bbb R^2$ with smooth boundary and
let $\widetilde \Gamma$ be an open set on $\partial\Omega$ and
$\Gamma_0 = \partial\Omega\setminus {\overline{\widetilde \Gamma}}$,
let $\nu$ be the unit outward normal vector to $\partial\Omega$.
 Consider the
following boundary value problem:
\begin{equation}\label{o-1}
L(x,D)u=\Delta u+2A\partial_z u+2B\partial_{\overline z} u+Qu=0\quad
\mbox{in}\,\,\Omega,\quad  u\vert_{\Gamma_0}=0.
\end{equation}
Here $u=(u_1,\dots, u_N)$ is an unknown vector-valued function
and $A, B, Q$ be smooth $N\times N$ matrices,
$i = \sqrt{-1}$, $x = (x_1, x_2) \in \Bbb R^2$, $x$ is identified with
$z = x_1 + ix_2 \in \Bbb C$,
$\partial_{z} = \frac{1}{2}\left(\frac{\partial}{\partial x_1}
- i\frac{\partial}{\partial x_2}\right)$ and
$\partial_{\overline z} = \frac{1}{2}\left(\frac{\partial}{\partial x_1}
+ i\frac{\partial}{\partial x_2}\right)$.

Consider the following partial Cauchy data:
$$
\mathcal C_{A,B,Q}=\{(u,\frac{\partial u}{\partial
\nu})\vert_{\widetilde \Gamma}; \thinspace
L(x,D)u=0\quad\mbox{in}\,\,\Omega,\quad u\vert_{\Gamma_0}=0,\,\,
u\in H^1(\Omega)\}.
$$

The paper is concerned with the following inverse problem: {\it
Using the partial Cauchy data $\mathcal C_{A,B,Q}$, determine matrix
coefficients $A,B,Q$. }

Note that we allowed freely choose Dirichlet data on
$\widetilde{\Gamma}$ and measure the corresponding $\frac{\partial
u}{\partial\nu}\vert_{\widetilde\Gamma}$. In one special case of
$N=1$ and $A=B=0$, this inverse boundary value problem is related to
so called the Calder\'on's problem (see \cite{C}), which is a
mathematical realization of {\sl Electrical Impedance Tomography.}

Similarly to the case of $N=1$ in \cite{IUY1}, the simultaneous
determination of all three coefficients $A,B,Q$ is impossible, but
we can establish some equations for coefficient matrices $(A,B,Q)$
which generate the same partial Cauchy data.

Our main result is
\begin{theorem}\label{vokal}
Let $A_j,B_j \in C^{5+\alpha}(\overline \Omega)$ and
$Q_j\in C^{4+\alpha}(\overline \Omega)$ for
$j=1,2$ and some $\alpha\in (0,1).$ Suppose that $\mathcal
C_{A_1,B_1,Q_1}=\mathcal C_{A_2,B_2,Q_2}.$ Then
\begin{equation}\label{op!}
A_1=A_2\quad\mbox{and}\,\, B_1=B_2\quad \mbox{on} \,\,\widetilde \Gamma,
\end{equation}
\begin{equation}\label{A1}
2\partial_z(A_1-A_2)+B_2(A_1-A_2)+(B_1-B_2)A_1-(Q_1-Q_2)=0\quad\mbox{in}
\,\,\Omega
\end{equation}
and
\begin{equation}\label{A2}
2\partial_{\overline
z}(B_1-B_2)+A_2(B_1-B_2)+(A_1-A_2)B_1-(Q_1-Q_2)=0\quad\mbox{in}\,\,\Omega.
\end{equation}
\end{theorem}

In the case of $N=1$ and two dimensions, there are many works and we
refer to some of them, and here we do not intend  to provide a
complete list. In the case $\widetilde \Gamma =\partial\Omega$ of
the full Cauchy data, the uniqueness in determining a potential $q$
in the two dimensional case was proved for the conductivity equation
by Nachman in \cite{N} within $C^4$ conductivities, and later in
\cite {AP} within $L^\infty$ conductivities.  For a convection
equation see \cite{ChengYama}.  The case of the Schr\"odinger
equation was solved by Bukhegim \cite {Bu}. In the case of the
partial Cauchy data on arbitrary subboundary, the uniqueness was
obtained in \cite{IUY} for potential $q \in
C^{5+\alpha}(\overline\Omega)$, and in \cite {IY}, the regularity
assumption was improved to $C^\alpha(\overline\Omega)$ in the case
of the full Cauchy data and up to $W^1_p(\Omega)$ with $p>2$ in the
case of partial Cauchy data on arbitrary subboundary. The case of
general second-order elliptic equation was studied in the papers
\cite{IUY1} and \cite{IUY2}. The results of \cite{IUY} were extended
to a Riemannian surface in \cite{GT}. The case where voltages are
applied and currents are measured on disjoint subboundaries was
discussed and the uniqueness is proved in \cite{IUY3}.   Conditional
stability estimates in determining a potential are obtained in
\cite{Nov}.  For the Calder\'on problem for the Schr\"odinger
equation in dimension three or more, we refer to the papers
\cite{BuU}, \cite{KSU}, \cite{KS} and \cite{SU}. To the best
knowledge of the authors, there are no publications for the
uniqueness for weakly coupling system of second-order elliptic
partial differential equations, and Theorem \ref{vokal} is the
affirmative answer.

Theorem \ref{vokal} asserts that any two coefficient matrices among three
are uniquely determined by partial Cauchy data on arbitrary
subboundary $\widetilde{\Gamma}$ for the system of elliptic differential
equations.
That is,
\begin{corollary}\label{coroA}
Let $(A_j,B_j,Q_j) \in C^{5+\alpha}(\overline\Omega) \times
C^{5+\alpha} (\overline\Omega)\times C^{4+\alpha}(\overline
\Omega)$, $j=1,2$ for some $\alpha\in (0,1)$ and be complex-valued.
We assume that either  $A_1\equiv A_2$ or  $B_1 \equiv B_2$ or $Q_1
\equiv Q_2$ in $\Omega$.  Then $\mathcal{C}_{A_1,B_1,Q_1}=
\mathcal{C} _{A_2,B_2,Q_2}$ implies $(A_1,B_1,Q_1) = (A_2,B_2,Q_2)$
in $\Omega$.
\end{corollary}

{\bf Proof.}

{\bf Case 1: $Q_1=Q_2.$}

Denote $R(x,D)(w_1,w_2)=(2\partial_z w_1+B_2w_1+w_2A_1,
2\partial_{\bar z}w_2+A_2w_2+w_1B_1).$ Therefore, applying
Theorem \ref{vokal}, we obtain
\begin{equation}\label{A10}
R(x,D)(A_1-A_2,B_1-B_2)=0\quad\mbox{in}\,\,\Omega
\end{equation}
and
\begin{equation}\label{A11}
(A_1-A_2)\vert_{\widetilde\Gamma}=(B_1-B_2)\vert_{\widetilde\Gamma}=0.
\end{equation}

Let a function $\psi\in C^2(\overline \Omega)$ satisfy $\vert
\nabla\psi\vert>0$ on $\overline \Omega$, $\lambda$ be a large
positive parameter and $\phi=e^{\lambda \psi}$. Then there exist
constants $\tau_0$ and $C$ independent of $\tau$ such that
\begin{equation}\label{voda1}
\vert \tau\vert^\frac 12 \Vert we^{\tau\phi}\Vert_{L^2(\Omega)} \le
C \Vert (\partial_z w)e^{\tau\phi}\Vert_{L^2(\Omega)},\quad \forall
\tau\ge \tau_0\,\,\mbox{and}\,\,\forall w\in H^1_0(\Omega)
\end{equation}
and
\begin{equation}\label{voda}
\vert \tau\vert^\frac 12 \Vert we^{\tau\phi}\Vert_{L^2(\Omega)} \le
C \Vert (\partial_{\overline z} w)e^{\tau\phi}\Vert_{L^2(\Omega)},
\quad \forall \tau\ge \tau_0 \,\,\mbox{and}\,\,\forall w\in
H^1_0(\Omega).
\end{equation}
Consider the boundary value problem
\begin{equation}\label{voda2}
R(x,D)(w_1,w_2)=(f_1,f_2)\quad \mbox{in}\,\,\Omega, \quad
(w_1,w_2)\vert_{\partial\Omega}=0.
\end{equation}
Applying the Carleman estimates (\ref{voda1}), (\ref{voda}) to each
of $N^2$ equations  in (\ref{voda2}), we have
\begin{equation}\label{dodo}
\vert \tau\vert^\frac 12 \Vert(
w_1,w_2)e^{\tau\phi}\Vert_{L^2(\Omega)} \le C
\left(\sum_{j=1}^2\Vert f_je^{\tau\phi}\Vert_{L^2(\Omega)} + \Vert
(w_1,w_2)e^{\tau\phi}\Vert_{L^2(\Omega)}\right), \quad \forall
\tau\ge \tau_0.
\end{equation}
The second term on the right-hand side of (\ref{dodo}) can be
absorbed into the left-hand side. Therefore we have
\begin{equation}\label{dodo1}
\vert \tau\vert^\frac 12 \Vert(
w_1,w_2)e^{\tau\phi}\Vert_{L^2(\Omega)} \le C \sum_{j=1}^2\Vert
f_je^{\tau\phi}\Vert_{L^2(\Omega)}, \quad \forall \tau\ge \tau_0.
\end{equation}
Using (\ref{dodo1}) and repeating the arguments in \cite{Ho}, we
prove that a solution of the Cauchy problem (\ref{A10}), (\ref{A11})
is zero.

{\bf Case 2: $B_1=B_2$.}

From equation (\ref{A2}), we have
$$
(A_1-A_2)B_1=(Q_1-Q_2) \quad\mbox{in}\,\,\Omega.
$$
Hence equation (\ref{A1}) can be written as
\begin{equation}\label{A13}
2\partial_z(A_1-A_2)+B_2(A_1-A_2)-(A_1-A_2)B_1=0\quad\mbox{in}\,\,\Omega,\quad
(A_1-A_2)\vert_{\widetilde\Gamma}=0.
\end{equation}
Using (\ref{voda1}), for the boundary value problem:
$$
2\partial_zw+B_2w-wB_1=f \quad\mbox{in}\,\,\Omega,
\quad w\vert_{\partial\Omega}=0,
$$
we obtain the estimate
\begin{equation}\label{dodo3}
\vert \tau\vert^\frac 12 \Vert we^{\tau\phi}\Vert_{L^2(\Omega)}
\le C \Vert f e^{\tau\phi}\Vert_{L^2(\Omega)} \quad \forall \tau\ge \tau_0.
\end{equation}
Using Carleman estimate (\ref{dodo3}) and repeating the arguments in
\cite{Ho}, we prove that solution of the Cauchy problem (\ref{A13}) is zero.
Then equation (\ref{A2}) implies that $Q_1=Q_2.$

The proof in the case $A_1=A_2$ is the same.$\blacksquare$

Next we consider other form of elliptic systems:
\begin{equation}\label{OM}
\widetilde L(x,D)u=\Delta u+\mathcal A\partial_{x_1} u +\mathcal
B\partial_{x_2} u + Qu.
\end{equation}
Here $\mathcal A$, $\mathcal B$, $Q$ are complex-valued $N\times N$
matrices. Let us define the following set of partial Cauchy data:
$$
\widetilde C_{\mathcal A,\mathcal B,Q}=\left\{(u, \frac{\partial
u}{\partial \nu})\vert_{\widetilde \Gamma}; \thinspace
\widetilde L(x,D)u=
\Delta u+\mathcal A\partial_{x_1} u +\mathcal B\partial_{x_2} u +
Qu=0\,\,\mbox{in}\,\,\Omega, \, u\vert_{\Gamma_0}=0, u\in
H^1(\Omega)\right\}.
$$
Then one can prove the following corollary.
\begin{corollary}\label{coroB}
Let $Q_1,Q_2\in C^{4+\alpha}(\overline\Omega)$ and let two pairs of
complex-valued coefficient matrices $(\mathcal A_{1}, \mathcal
B_{1}), (\mathcal A_{2}, \mathcal B_{2}) \in C^{5+\alpha}(\overline
\Omega)\times C^{5+\alpha}(\overline \Omega)$  for some  $\alpha\in
(0,1).$
We assume that $Q_1\equiv Q_2$ in $\Omega$. Then
$(\mathcal A_{1}, \mathcal B_{1})\equiv (\mathcal A_{2}, \mathcal
B_{2})$ in $\Omega$.
\end{corollary}

{\bf Proof.}  Observe that $\widetilde L(x,D)=\Delta + A\partial_{z}
+B\partial_{\bar z}  + Q$ where $A=\mathcal A+i\mathcal B$ and
$B=\mathcal A-i\mathcal B.$ Therefore, applying Corollary 2, we
complete the proof. $\blacksquare$

{\bf Remark. } Unlike Corollary \ref{coroA}, in the two cases of
$\mathcal A_1\equiv\mathcal A_2$ and $\mathcal B_1\equiv\mathcal B_2$,
we can not, in general, claim that
$(\mathcal A_1,\mathcal B_1, Q_1)= (\mathcal A_2,\mathcal B_2, Q_2).$
By the same argument as Corollary \ref{coroA}, we can prove only \\
(i) $\frac{\partial \mathcal B_1}{\partial x_1} =
\frac{\partial \mathcal B_2}{\partial x_1}$ in $\Omega$ if
$\mathcal A_1 = \mathcal A_2$ in $\Omega$.
\\
(ii) $\frac{\partial \mathcal A_1}{\partial x_2} =
\frac{\partial \mathcal A_2}{\partial x_2}$ in $\Omega$
if $\mathcal B_1 = \mathcal B_2$ in $\Omega$.

Moreover consider the following example
$$
\Omega = (0,1)\times (0,1),
$$
$$
\widetilde \Gamma=\{(x_1,x_2); \thinspace x_2=0, \,\,  0<x_1<1\}\cup
\{(x_1,x_2)\vert x_2=1, \,\,  0<x_1<1\},
$$
and let us choose $\eta(x_2)\in C^\infty_0(0,1)$.
Then the operators $\widetilde L(x,D)$ and $e^{s\eta}\widetilde
L(x,D)e^{-s\eta}$
generate the same partial Cauchy data, but the matrix coefficient matrices
are not equal.

\section{Preliminary results}\label{sec1}

Throughout the paper, we use the following notations.
\\

\noindent {\bf Notations.} $i=\sqrt{-1}$, $x_1, x_2, \xi_1, \xi_2
\in {\Bbb R}^1$, $z=x_1+ix_2$, $\zeta=\xi_1+i\xi_2$, $\overline{z}$
denotes the complex conjugate of $z \in \Bbb C$. We identify $x =
(x_1,x_2) \in {\Bbb R}^2$ with $z = x_1 +ix_2 \in {\Bbb C},$
$\partial_z = \frac 12(\partial_{x_1}-i\partial_{x_2})$,
$\partial_{\overline z}= \frac12(\partial_{x_1}+i\partial_{x_2}),$
$\beta=(\beta_1,\beta_2), \vert \beta\vert=\beta_1+\beta_2.$
$D=(\frac{1}{i}\frac{\partial}{\partial
x_1},\frac{1}{i}\frac{\partial}{\partial x_2}).$ Let $\chi_G$ be the
characteristic function of  the set $G.$ The tangential derivative
on the boundary is given by
$\partial_{\vec\tau}=\nu_2\frac{\partial}{\partial x_1}
-\nu_1\frac{\partial}{\partial x_2},$ where $\nu=(\nu_1, \nu_2)$ is
the unit outer normal to $\partial\Omega,$ $B(\widehat
x,\delta)=\{x\in \Bbb R^2; \thinspace \vert x-\widehat x\vert<
\delta\},$ $S(\widehat x,\delta)=\{x\in \Bbb R^2; \thinspace \vert
x-\widehat x\vert= \delta\}$. We set $(u,v)_{L^2(\Omega)} =
\int_{\Omega} u\overline{v} dx$ for functions $u, v$, while by
$(a,b)$ we denote the scalar product in $\Bbb R^2$ if there is no
fear of confusion. For $f:{\Bbb R}^2\rightarrow {\Bbb R}^1$, the
symbol $f''$ denotes the Hessian matrix  with entries
$\frac{\partial^2 f}{\partial x_k\partial x_j},$ $\mathcal L(X,Y)$
denotes the Banach space of all bounded linear operators from a
Banach space $X$ to another Banach space $Y$.  Let $E$ be the
$N\times N$ unit matrix. We set $ \Vert u\Vert_{H^{1,\tau}(\Omega)}
= (\Vert u\Vert_{H^1(\Omega)}^2 + \vert \tau\vert^2\Vert
u\Vert^2_{L^2(\Omega)})^{\frac{1}{2}}$. Finally for any $\widetilde
x\in \partial\Omega$, we introduce the left and the right tangential
derivatives as follows:
$$
{\mathbf D}_+(\widetilde x) f=\lim_{s\rightarrow +0}
\frac{f(\ell(s))-f(\widetilde x)}{s},
$$
where $\ell(0)=\widetilde x,$ $\ell(s)$ is a parametrization of
$\partial\Omega$ near $\widetilde x$, $s$ is the length of the
curve, and we are moving clockwise as $s$ increases;
$$
{\mathbf D}_-(\widetilde x) f=\lim_{s\rightarrow -0}
\frac{f(\widetilde \ell(s))-f(\widetilde x)}{s},
$$
where $\widetilde\ell(0)=\widetilde x,$ $\widetilde \ell(s)$ is the
parametrization of $\partial\Omega$ near $\widetilde x$ , $s$ is the
length of the curve, and we are moving counterclockwise as $s$
increases. By $o_X(\frac{1}{\tau^\kappa})$ we denote a function
$f(\tau,\cdot)$ such that $ \Vert
f(\tau,\cdot)\Vert_X=o(\frac{1}{\tau^\kappa})\quad \mbox{as}
\,\,\vert \tau\vert\rightarrow +\infty. $\vspace{0.4cm}

For some $\alpha\in (0,1)$, we consider a function
$\Phi(z)=\varphi(x_1,x_2)+i\psi(x_1,x_2) \in
C^{6+\alpha}(\overline{\Omega})$ with real-valued $\varphi$ and
$\psi$ such that
\begin{equation}\label{zzz}
\partial_z\Phi(z) = 0 \quad \mbox{in}
\,\,\Omega, \quad\mbox{Im}\,\Phi\vert_{\Gamma_0^*}=0,
\end{equation}
where $\Gamma_0^*$ is an open set on $\partial\Omega$ such that
$\Gamma_0\subset\subset \Gamma_0^*.$ Denote by $\mathcal H$ the set
of all the critical points of the function $\Phi$:
$$
\mathcal H = \{z\in\overline\Omega; \thinspace
\frac{\partial\Phi}{\partial z} (z)=0\}.
$$
Assume that $\Phi$ has no critical points on
$\overline{\widetilde\Gamma}$, and that all critical points  are
nondegenerate:
\begin{equation}\label{mika}
\mathcal H\cap \partial\Omega\subset\Gamma_0,\quad
\partial^2_{z}\Phi (z)\ne 0, \quad \forall z\in
\mathcal H.
\end{equation}
Then $\Phi$  has only a finite number of critical points and we can
set:
\begin{equation}\label{mona}
{\mathcal H} \setminus \Gamma_0= \{ \widetilde{x}_1, ...,
\widetilde{x}_{\ell} \},\quad \mathcal H \cap \Gamma_0=\{
\widetilde{x}_{\ell+1}, ..., \widetilde{x}_{\ell+\ell'} \}.
\end{equation}

Let $\partial \Omega=\cup_{j=1}^{\mathcal N}\gamma_j,$ where
$\gamma_j$ is a closed contour. The following proposition was proved
in \cite{IUY}.

\begin{proposition}\label{Proposition -1}
Let $\widetilde x$ be an arbitrary point in $\Omega.$ There exists a
sequence of functions $\{\Phi_\epsilon\}_{\epsilon\in(0,1)}$
satisfying (\ref{zzz}), (\ref{mika}) and there exists a sequence
$\{\widetilde x_\epsilon\}, \epsilon\in (0,1)$ such that
$$
\widetilde x_\epsilon \in \mathcal H_\epsilon
= \{z\in\overline\Omega; \thinspace
\frac{\partial \Phi_\epsilon}{\partial
z}(z)=0 \},\quad \widetilde x_\epsilon\rightarrow \widetilde
x\,\,\mbox{ as}\,\, \epsilon\rightarrow +0.
$$
Moreover for any $j$ from $\{1,\dots,\mathcal N\}$, we have
$$
\mathcal H_\epsilon\cap\gamma_j=\emptyset\quad\mbox{if}\,\,
\gamma_j\cap \widetilde \Gamma\ne\emptyset,
$$
$$
\mathcal H_\epsilon\cap\gamma_j\subset \Gamma_0 \quad\mbox{if}\,\,
\gamma_j\cap \widetilde \Gamma = \emptyset
$$
and
$$
\mbox{Im}\,\Phi_\epsilon(\widetilde x_\epsilon)\notin \{\mbox{Im}\,
\Phi_\epsilon(x); \thinspace x\in \mathcal H_\epsilon\setminus
\{\widetilde{x_\epsilon}\}\} \,\,\mbox{and}
\,\,\mbox{Im}\,\Phi_\epsilon(\widetilde x_\epsilon) \ne 0.
$$
\end{proposition}

The following proposition  was proved in \cite{IUY1}.
\begin{proposition}\label{Proposition -2}
Let $\widehat\Gamma_*\subset\subset\widetilde \Gamma$ be an arc with
the left endpoint $x_-$ and the right endpoint $x_+$ oriented clockwise.
For any $\widehat x\in Int\,\widehat\Gamma_*$, there exists a function
$\Phi(z)$ which satisfies (\ref{zzz}), (\ref{mika}),
$\mbox{Im}\,\Phi\vert_{\partial\Omega\setminus {\widehat\Gamma_{*}}}=0$,
\begin{equation}\label{lana0}
\widehat x\in\mathcal G=\{x\in\widehat\Gamma_*; \quad
\frac{\partial  \mbox{Im}\,\Phi}{\partial \vec\tau}( x)=0\}, \quad
card\, \mathcal G <\infty
\end{equation}
and
\begin{equation}
(\frac{\partial}{\partial \vec\tau})^2 \mbox{Im}\,\Phi(x)\ne 0 \quad
\forall x\in \mathcal G\setminus\{x_-,x_+\}.
\end{equation}
Moreover
\begin{equation}\label{lana1}
\mbox{Im}\,\Phi(\widehat x)\ne \mbox{Im}\,\Phi(x),\quad \forall x\in
\mathcal G\setminus \{\widehat x\} \mbox { and }\quad
\mbox{Im}\,\Phi(\widehat x)\neq 0
\end{equation}
and
\begin{equation}\label{zopa}{\mathbf D}_+(x_-)
(\frac{\partial}{\partial \vec \tau})^6\mbox{Im}\,\Phi\ne 0, \quad
{\mathbf D}_-(x_+)(\frac{\partial}{\partial \vec \tau})^6
\mbox{Im}\,\Phi\ne 0.
\end{equation}
\end{proposition}

Later we use the following Proposition (see \cite{IUY}) :

\begin{proposition}\label{gandonnal} Let $\Phi$ satisfy (\ref{zzz})
and (\ref{mika}). For every $g\in L^1(\Omega)$, we have
$$
\int_\Omega ge^{\tau(\Phi-\overline \Phi)}dx\rightarrow \,\,0\quad
\mbox{as}\quad \tau\rightarrow +\infty .
$$
\end{proposition}

Moreover
\begin{proposition}\label{gandon} Let $\Phi$  satisfy (\ref{zzz}),
(\ref{mika}), $g\in W^1_p(\Omega)$ with some $p> 2$,
$g\vert_{\mathcal H}=0$ and $supp \,g\subset \Omega.$ Then
$$
\int_\Omega ge^{\tau(\Phi-\overline \Phi)}dx=o(\frac 1\tau)\quad
\mbox{as}\quad \tau\rightarrow +\infty .
$$
\end{proposition}

{\bf Proof.} By the Sobolev imbedding theorem, the function $g$
belongs to $C^\alpha(\overline \Omega)$ for some positive $\alpha.$ Note
that by (\ref{mika}) and the assumption on $g$, we have
\begin{equation}\label{mmm} \Vert g\frac{(\nabla \psi, \nu)
e^{\tau(\Phi-\overline
\Phi)}}{2i\vert\nabla \psi\vert^2}\Vert_{C(\overline{S(\widetilde
x_j,\delta)})}\le \frac{C\Vert g\Vert
_{C(\overline{S(\widetilde x_j,\delta)})}}{\delta}
\le \frac{C}{\delta^{1-\alpha}}.
\end{equation}
Also
$$
\mbox{div}\,(g\frac{\nabla \psi}{2i\tau\vert\nabla \psi\vert^2})
= (\nabla g,\frac{\nabla \psi}{2i\tau\vert\nabla
\psi\vert^2})
+ g \mbox{div}\left(\frac{\nabla \psi}{2i\tau\vert\nabla
\psi\vert^2}\right).
$$
Since
$$
\vert (\nabla g,\frac{\nabla
\psi}{2i\tau\vert\nabla \psi\vert^2})\vert \le C
\sum_{j=1}^{\ell+\ell'}\frac{\vert \nabla g (x)\vert}{\vert x-\widetilde
x_j\vert},
$$
by the H\"older inequality we conclude that $(\nabla
g,\frac{\nabla \psi}{2i\tau\vert\nabla \psi\vert^2})\in L^1(\Omega)$.
By (\ref{mika}) and assumption that $g\vert_{\mathcal H}=0$, we obtain
$$
\left\vert g \mbox{div}\left(\frac{\nabla
\psi}{2i\tau\vert\nabla \psi\vert^2}\right)\right\vert
\le C \sum_{j=1}^{\ell+\ell'}\frac{\vert  g(x)\vert}
{\vert x-\widetilde x_j\vert^2}
\le C\sum_{j=1}^{\ell+\ell'}\frac{\Vert
g\Vert_{C^\alpha(\overline \Omega)} }{\vert x-\widetilde
x_j\vert^{2-\alpha}}.
$$
Therefore $\mbox{div}\left(g\frac{\nabla
\psi}{2i\tau\vert\nabla \psi\vert^2}\right)\in L^1(\Omega)$.
By (\ref{mmm}), passing to the limit as $\delta$ goes to zero, we have
$$
J=\int_\Omega ge^{\tau(\Phi-\overline \Phi)}dx=\lim_{\delta\rightarrow
0}\int_{\Omega\setminus \cup_{j=1}^{\ell+\ell'} S(\widetilde
x_j,\delta)}ge^{\tau(\Phi-\overline \Phi)}dx=\lim_{\delta\rightarrow
0}\int_{\Omega\setminus \cup_{j=1}^{\ell+\ell'} S(\widetilde
x_j,\delta)}g\frac{(\nabla \psi, \nabla)e^{\tau(\Phi-\overline
\Phi)}}{2i\tau\vert\nabla \psi\vert^2}dx
$$
$$
= \lim_{\delta\rightarrow 0}\int_{\cup_{j=1}^{\ell+\ell'}
S(\widetilde x_j,\delta)}g\frac{(\nabla \psi, \nu)
e^{\tau(\Phi-\overline
\Phi)}}{2i\tau\vert\nabla
\psi\vert^2}d\sigma-\lim_{\delta\rightarrow 0}\int_{\Omega\setminus
\cup_{j=1}^{\ell+\ell'} S(\widetilde
x_j,\delta)}\mbox{div}\,(g\frac{\nabla \psi}{2i\tau\vert\nabla
\psi\vert^2})e^{\tau(\Phi-\overline
\Phi)}dx
$$
$$
= -\int_{\Omega}\mbox{div}\,(g\frac{\nabla
\psi}{2i\tau\vert\nabla \psi\vert^2})e^{\tau(\Phi-\overline \Phi)}dx.
$$
Using Proposition \ref{gandonnal}, we finish the
proof. $\blacksquare$

Consider the boundary value problem
$$
L(x,{D})u = f \quad \mbox{in} \quad \Omega, \quad u \vert_{\partial
\Omega} = 0.
$$
\begin{proposition}\label{Theorem 2.1}
Suppose that $\Phi$ satisfies (\ref{zzz}), (\ref{mika}), $u\in
H^1_0(\Omega)$ and $\Vert A\Vert _{L^\infty(\Omega)}+\Vert
B\Vert _{L^\infty(\Omega)}+\Vert
Q\Vert_{L^\infty(\Omega)}\le K$. Then there exist
$\tau_0=\tau_0(K,\Phi)$ and $C=C(K,\Phi)$, independent of $u$ and
$\tau$, such that
\begin{eqnarray}\label{suno4} \vert \tau\vert\Vert
ue^{\tau\varphi}\Vert^2_{L^2(\Omega)}+\Vert
ue^{\tau\varphi}\Vert^2_{H^1(\Omega)}+\Vert\frac{\partial
u}{\partial\nu}e^{\tau\varphi}\Vert^2_{L^2(\Gamma_0)}+
\tau^2\Vert\vert\frac{\partial\Phi}{\partial z} \vert
ue^{\tau\varphi}\Vert^2_{L^2(\Omega)}\nonumber \\
\le C(\Vert (L(x,D)
u)e^{\tau\varphi}\Vert^2_{L^2(\Omega)}+\vert\tau\vert
\int_{\widetilde\Gamma^*}\vert \frac{\partial
u}{\partial\nu}\vert^2e^{2\tau\varphi}d\sigma)
\end{eqnarray}
for all $\vert\tau\vert>\tau_0$.
\end{proposition}
For the scalar equation, the estimate is proved in \cite{IUY1} .
In order to prove this estimate for the system, it is sufficient to
apply the scalar estimate to each equation in the system and take an
advantage of the second large parameter in order to absorb the
right-hand side.

Using estimate (\ref{suno4}), we obtain
\begin{proposition} \label{vanka}
There exists a constant $\tau_0$ such that for $\vert
\tau\vert\ge \tau_0$ and any $f\in L^2(\Omega)$, there exists a
solution to the boundary value problem
\begin{equation}\label{lola}
L(x,D)u =f\quad\mbox{in}\,\,\Omega, \quad u\vert_{\Gamma_0}=0
\end{equation}
such that
\begin{equation}\label{2}
\Vert u\Vert_{H^{1,\tau}(\Omega)}/\root\of{\vert\tau\vert}
\le C\Vert  f\Vert_{L^2(\Omega)}.
\end{equation}
Moreover if $f/\partial_z\Phi\in L^2(\Omega)$, then for any $\vert
\tau\vert\ge \tau_0$  there exists a solution to the boundary value
problem (\ref{lola}) such that
\begin{equation}\label{3}
\Vert u\Vert_{H^{1,\tau}(\Omega)} \le C\Vert
f/\partial_z\Phi\Vert_{L^2(\Omega)}.
\end{equation}
The constants $C$ in (\ref{2}) and (\ref{3}) are independent of
$\tau.$
\end{proposition}

The proof is exactly the same as the proof of Proposition 2.5  in
\cite{IUY1} and relies on  the Carleman estimate (\ref{suno4}).

Let us introduce the operators:
$$
\partial_{\overline z}^{-1}g=-\frac 1\pi\int_\Omega
\frac{g(\xi_1,\xi_2)}{\zeta-z}d\xi_1d\xi_2,\quad
\partial_{ z}^{-1}g=-\frac 1\pi\int_\Omega
\frac{g(\xi_1,\xi_2)}{\overline\zeta-\overline z}d\xi_1d\xi_2.
$$

Then we have (e.g., p.47, 56, 72 in \cite{VE}):
\begin{proposition}\label{Proposition 3.0}
{\bf A)} Let $m\ge 0$ be an integer number and $\alpha\in (0,1).$
Then $\partial_{\overline z}^{-1},\partial_{ z}^{-1}\in \mathcal
L(C^{m+\alpha}(\overline{\Omega}),C^{m+\alpha+1}
(\overline{\Omega})).$
\newline
{\bf B}) Let $1\le p\le 2$ and $ 1<\gamma<\frac{2p}{2-p}.$ Then
 $\partial_{\overline z}^{-1},\partial_{ z}^{-1}\in
\mathcal L(L^p( \Omega),L^\gamma(\Omega)).$
\newline
{\bf C})Let $1< p<\infty.$ Then  $\partial_{\overline z}^{-1},
\partial_{ z}^{-1}\in \mathcal L(L^p( \Omega),W^1_p(\Omega)).$
\end{proposition}

For any matrix $B\in C^{5+\alpha}(\overline \Omega)$, consider the linear
operators $T_B$ and $P_B$ such that
\begin{equation}
(2\partial_z+B)T_Bg=g\quad \mbox{in}\,\,\Omega; \quad
(2\partial_{\overline z}+B)P_Bg=g\quad \mbox{in}\,\,\Omega
\end{equation}
and
\begin{equation}\label{sonka}
T_B,P_B\in \mathcal L(H^s(\Omega), H^{s+1}(\Omega))\cap \mathcal
L(C^{k+\alpha}(\Omega), C^{k+1+\alpha}(\Omega))\quad \forall s\in
[0,6], \forall k\in \{0,1,\dots 6\},
\end{equation}

and
\begin{equation}\label{sonka1}
T_B,P_B\in \mathcal L((H^1(\Omega))', L^2(\Omega)).
\end{equation}
 The existence of the operators $T_B, P_B$ with the above
properties follows from the regularity theory of elliptic systems on
the plane (see e.g., \cite{Wendland}).

Let $e\in C^\infty_0(\Omega)$ satisfy $\vert e(x)\vert\le 1$,
the support of $e$ be concentrated in a small neighborhood of
$\mathcal H\setminus \overline \Gamma_0$ and $e$ be identically equal to one
in an open set $\mathcal O$ which contains $\mathcal H\setminus
\overline \Gamma_0.$
We introduce the operators $\frak T_B$ and $ \frak P_B$ by
\begin{equation}\label{o11}
\frak T_B=\frac 12\sum_{j=0}^\infty (-1)^j(\frac 12 \partial_z^{-1}
e B)^j
\partial_z^{-1}, \quad \frak P_B
= \frac 12\sum_{j=0}^\infty (-1)^j(\frac 12
\partial_{\overline z}^{-1} e B)^j \partial_{\overline z}^{-1}.
\end{equation}
Taking the function $e$ such that $\int_{\mbox{supp}\, e}1dx$ is sufficiently
small, we have
\begin{equation}\label{tot}
\Vert \partial_z^{-1} e B\Vert_{\mathcal L(L^p(\Omega),
L^p(\Omega))}<1\quad \mbox{and}\quad \Vert \partial_{\overline z}^{-1} e
B\Vert_{\mathcal L(L^p(\Omega), L^p(\Omega))}<1.
\end{equation}
Indeed, by Proposition \ref{Proposition 3.0} for any $p>1$
there exists a number $q\in (1,p)$ such that the operators
$\partial_z^{-1}, \partial_{\overline z}^{-1}:L^q(\Omega)
\rightarrow L^p(\Omega)$ are continuous.  Therefore
$$
\Vert \partial_z^{-1}eBg\Vert_{L^p(\Omega)}
\le \Vert \partial_z^{-1}\Vert_{\mathcal L(L^q(\Omega),L^p(\Omega))}
\Vert B\Vert_{L^\infty(\Omega)}\Vert eg\Vert_{L^q(\Omega)}
$$
$$
\le \Vert \partial_z^{-1}\Vert_{\mathcal L(L^q(\Omega),L^p(\Omega))}
\Vert B\Vert_{L^\infty(\Omega)}\left(\int_{supp\, e}1dx\right)^{(p-q)/p}
\Vert g\Vert_{L^p(\Omega)},
$$
and if $\int_{\mbox{supp}\, e}1dx$ is small, then we easily
have (\ref{tot}).

Hence the operators $\frak T_B$ and $\frak P_B$ introduced in
(\ref{o11}) are correctly defined.

We define two other operators:
\begin{equation}\label{anna}
\mathcal R_{\tau}g = \frac 12e^{\tau(\Phi-\overline {\Phi})}
\partial_{\overline z}^{-1}(g
e^{\tau(\overline{\Phi}- {\Phi})}),\,\, \widetilde {\mathcal
R}_{\tau}g = \frac 12e^{\tau(\overline {\Phi}-{\Phi})}
\partial_{ z}^{-1}(g e^{\tau( {\Phi}
-\overline {\Phi})}).
\end{equation}
For any $N\times N$ matrix $B$ with elements from $C^1(\overline\Omega)$,
we set
$$
\mbox{T}_B=\frak T_{B}-T_B(1-e)B\frak T_{B},\quad   \mbox{P}_B=\frak
P_{B}-P_B(1-e)B\frak
P_{B},
$$
\begin{equation}
\widetilde {\mathcal R}_{\tau, B}g=\frak T_{B,\tau}
g-e^{\tau(\overline\Phi-\Phi)}T_B(e^{\tau(\Phi-\overline\Phi)}(1-e)B\frak
T_{B,\tau} g),
\end{equation}
$$
{\mathcal R}_{\tau, B}g=\frak P_{B,\tau}
g-e^{\tau(\Phi-\overline\Phi)}P_B(e^{\tau(\overline\Phi-\Phi)}(1-e) B\frak
P_{B,\tau} g )
$$
and
\begin{equation}
\frak T_{B,\tau}=e^{\tau(\overline\Phi-\Phi)}\frak
T_{B}e^{\tau(\Phi-\overline\Phi)},\quad \frak
P_{B,\tau}=e^{\tau(\Phi-\overline\Phi)}\frak
P_{B}e^{\tau(\overline\Phi-\Phi)}.
\end{equation}

For any $g \in C^\alpha(\overline\Omega)$, the functions
$\mathcal R_{\tau, B}g$ and $\widetilde {\mathcal R}_{\tau, B}g$
solve the equations:
\begin{equation}\label{0011}
(2\partial_{\overline z}+2\tau \partial_{\overline z}\overline \Phi+B) \mathcal
R_{\tau, B} g=g\quad \mbox{in}\,\,\Omega,\quad (2\partial_{ z}+2\tau
\partial_{ z} \Phi+B)\widetilde {\mathcal R}_{\tau, B} g=g\quad
\mbox{in}\,\,\Omega .
\end{equation}

We have
\begin{proposition}\label{vasya}
Let $B\in C^1(\overline\Omega), g\in C^2(\overline \Omega),
\mbox{supp}\, g\subset\subset \{x\vert e(x)=1\}$ and $g\vert_{\mathcal H}=0$.
Then for $p\in (1,\infty)$, we have
\begin{equation}
\Vert \widetilde {\mathcal R}_{\tau,
B}g-\frac{g}{2\tau\partial_z\Phi}\Vert_{L^p(\Omega)}+\Vert {\mathcal
R}_{\tau, B}g-\frac{g}{2\tau\partial_{\overline
z}\overline\Phi}\Vert_{L^p(\Omega)}=o(\frac 1\tau)\quad
\mbox{as}\,\,\vert\tau\vert\rightarrow \infty.
\end{equation}
\end{proposition}
{\bf Proof.}
By Proposition 3.4 of \cite{IUY1}, for any $p>1$, we have
\begin{equation}\label{002}
\Vert \widetilde {\mathcal
R}_{\tau}g-\frac{g}{2\tau\partial_z\Phi}\Vert_{L^p(\Omega)}+\Vert
{\mathcal R}_{\tau}g-\frac{g}{2\tau\partial_{\overline
z}\overline\Phi}\Vert_{L^p(\Omega)}=o(\frac 1\tau)\quad
\mbox{as}\,\,\vert\tau\vert\rightarrow \infty.
\end{equation}
Propositions \ref{gandon} and \ref{Proposition 3.0} yield
\begin{equation}\label{001}
\Vert R_{\tau}\{\frac{g}{2\tau\partial_z\Phi}\}\Vert
_{L^p(\Omega)}
+ \Vert\widetilde
R_{\tau}\{\frac{g}{2\tau\partial_{\overline
z}\overline\Phi}\}\Vert_{L^p(\Omega)}=o(\frac{1}{\tau})\quad
\mbox{as}\,\,\vert\tau\vert\rightarrow \infty.
\end{equation}

Thanks to (\ref{001}) and (\ref{002}), we obtain
\begin{equation}\label{popo!}
\Vert \frak T_{B,\tau} g-\frac{g}{2\tau\partial_z\Phi}\Vert_{L^p(\Omega)}
+ \Vert \frak P_{B,\tau} g-\frac{g}{2\tau\partial_{\overline
z}\overline\Phi}\Vert_{L^p(\Omega)}=o(\frac 1\tau)\quad
\mbox{as}\,\,\vert\tau\vert\rightarrow \infty.
\end{equation}
By $\mbox{supp}\,g\subset\subset \{x\vert e(x)=1\}$ and
(\ref{sonka}), (\ref{popo!}), we obtain the asymptotic formula:
$$
\Vert e^{\tau(\overline \Phi-\Phi )}T_B e^{\tau(\Phi-\overline\Phi)}\circ
(1-e)B \frak T_{B,\tau}
g \Vert_{L^p(\Omega)}+\Vert e^{\tau( \Phi-\overline\Phi )}P_B
e^{\tau(\overline\Phi-\Phi)}\circ (1-e)B
 \frak P_{B,\tau}
 g \Vert_{L^p(\Omega)}=o(\frac {1}{\tau})
\quad \mbox{as}\,\,\vert\tau\vert \rightarrow \infty.
$$
The proof is completed. $\blacksquare$
%
%
%
%
%
%
\section{Proof of Theorem \ref{vokal}}\label{sec2}

{\bf Proof  of Theorem \ref{vokal}.}

{\bf Step 1: Construction of complex geometric optics solutions.}

Let the function $\Phi$ satisfy (\ref{zzz}), (\ref{mika}) and $\widetilde x$ be
some point from $\mathcal H\setminus\Gamma_0.$
Without loss of generality, we may assume that $\widetilde \Gamma$ is an arc
with the endpoints $x_\pm.$

Consider the following operator:
\begin{eqnarray}\label{ooo} L_1(x,D)=4\partial_z\partial_{\overline z}+2
A_1\partial_z+2B_1\partial_{\overline z}+Q_1 \nonumber\\
= (2\partial_z+B_1)(2\partial_{\overline z}+A_1)+Q_1-2\partial_z A_1-B_1A_1
\nonumber\\
= (2\partial_{\overline z}+A_1)(2\partial_z+B_1)+Q_1-2\partial_{\overline z}
B_1-A_1B_1 .
\end{eqnarray}

Let $(w_{0},\widetilde w_{0})\in C^{6+\alpha}(\overline \Omega)$
be a nontrivial solution to the boundary value problem:
\begin{equation}\label{-5}
\mathcal K(x,D)(w_{0},\widetilde w_{0})=(2\partial_{\overline z}w_{0} +A_1
w_{0}, 2\partial_{ z}\widetilde w_{0} +B_1 \widetilde
w_{0})=0\quad\mbox{in}\,\,\Omega, \quad w_{0}+\widetilde w_{0}=0\quad
\mbox{on}\,\,\Gamma_0.
\end{equation}

We have
\begin{proposition}\label{nikita}
Let $\widetilde x$ be an arbitrary
point from $\mathcal H\setminus \overline\Gamma_0$ and
$\vec z\in \Bbb C^N$ be an arbitrary vector.
There exists a solution $(w_{0},\widetilde
w_{0})\in C^{6+\alpha}(\overline \Omega)$ to problem (\ref{-5}) such that
\begin{equation}\label{xoxo1}
w_{0}(\widetilde x)=\vec z,
\end{equation}
\begin{equation}\label{xoxo1u}
\lim_{x\rightarrow x_\pm}\frac{ \vert w_0(x)\vert}{\vert x-x_\pm\vert^{98}}
=\lim_{x\rightarrow x_\pm}\frac{ \vert \widetilde w_0(x)\vert}
{\vert x-x_\pm\vert^{98}}=0
\end{equation}
and
\begin{equation}\label{xoxo}
\partial^{\alpha_1}_{x_1}\partial^{\alpha_2}_{x_2} w_{0}(x)
= \partial^{\alpha_1}_{x_1}\partial^{\alpha_2}_{x_2} \widetilde w_{0}(x)
\quad \forall x\in
\mathcal H\setminus\{\widetilde x\}\,\,\, \mbox{and}\,\,\, \forall
\alpha_1+\alpha_2\le 6.
\end{equation}
\end{proposition}

{\bf Proof.} Let us fix a point $\widetilde x$ from $\mathcal
H\setminus \{\widetilde x\} .$ By Proposition 4.2 of \cite{IUY1}
there exists a holomorphic function $a(z)\in C^7(\overline\Omega)$
such that $Im\, a\vert_{\Gamma_0}=0,$ $a(\widetilde x)=1$ and $a$
vanishes at each point of the set $\{x_\pm\}\cup\mathcal H\setminus
\{\widetilde x\}$. Let $(w_{0,0} ,\widetilde w_{0,0})\in
C^{6+\alpha}(\overline\Omega)$ be a solution to problem (\ref{-5})
such that $w_{0,0}(\widetilde x)=\vec z.$ Since $(w_0, \widetilde
w_0)=(a^{10}w_{0,0},\overline a^{10} \widetilde w_{0,0})$ solves
problem (\ref{-5}) and satisfies (\ref{xoxo}) -(\ref{xoxo1}), the
proof of the proposition is completed. $\blacksquare$

Now we start the construction of complex geometric optics solution.
Let the pair $(w_0,\widetilde w_0)$ be defined by Proposition
\ref{nikita}.  Short computations and (\ref{ooo})  yield
\begin{equation}\label{oi}
L_1(x,D) (w_0e^{\tau\Phi})=(Q_1-2\partial_z
A_1-B_1A_1)w_0e^{\tau\Phi},\quad L_1(x,D) (\widetilde
w_0e^{\tau\overline \Phi})=(Q_1-2\partial_{\overline z}
B_1-A_1B_1)\widetilde w_0 e^{\tau\overline\Phi}. \end{equation}

Let $e_1,e_2$ be smooth functions such that
\begin{equation}\label{short}
\mbox{supp}\,e_1\subset\subset\mbox{supp}\,e=1,\quad e_1+e_2=1\quad
\mbox{on}\,\,\Omega,\quad
\end{equation}
and $e_1$ vanishes in a neighborhood of $\partial\Omega$ and $e_2$
vanishes in a neighborhood of the set $\mathcal H\setminus \overline
\Gamma_0.$

Denote $G_\epsilon=\{x\in \Omega\vert dist (\mbox{supp}\,
e_1,x)>\epsilon\}.$ We have
\begin{proposition} \label{popo}
Let $B, q\in C^{5+\alpha}(\overline\Omega)$  for some positive
$\alpha$  and $\widetilde q\in W^1_p(\overline\Omega)$ for some
$p>2.$ Suppose that $q\vert_{\mathcal H}=\widetilde q\vert_{\mathcal
H}=0.$  There exist smooth functions $m_\pm\in C^2(\overline
G_\epsilon)$ which is independent of $\tau$ such that for any
$\overline G_\epsilon\cap \mbox{supp}\,e =\emptyset$, the asymptotic
formulae hold true:
\begin{eqnarray}\label{50}
\widetilde{\mathcal R}_{\tau, B_1}(e_1(q+\frac{\widetilde q}{\tau}))
= e^{\tau(\overline\Phi-\Phi)}\left (\frac{m_+ e^{2i\tau\psi
(\widetilde x)}}{\tau^2}+o_{C^2(\overline
G_\epsilon)}(\frac{1}{\tau^2})\right)
\quad\mbox{as}\,\vert\tau\vert\rightarrow +\infty ,\\
\quad{\mathcal R}_{\tau, A_1} (e_1(q + \frac{\widetilde q}{\tau}))
=e^{\tau(\Phi-\overline\Phi)}\left (\frac{m_-
e^{-2i\tau\psi(\widetilde x)}}{\tau^2} +o_{C^2(\overline
G_\epsilon)}(\frac{1}{\tau^2})\right
)\quad\mbox{as}\,\vert\tau\vert\rightarrow +\infty.
\end{eqnarray}
\end{proposition}

{\bf Proof.} By the Sobolev imbedding theorem the function $\tilde
q$ belong to the space $C^\alpha(\bar\Omega)$ with some positive
$\alpha.$
 Therefore the trace of $\tilde q$ on $\mathcal H$ defined correctly. For all $N$ and for any domain
$G_{\epsilon_0} $  with $\epsilon_0 >0$ there exists a function
$m_{+,N}\in C^2(\overline G_{\epsilon_0})$ such that
\begin{equation}\label{ika}
\frac {e^{-2\tau i\psi}}{2} (-1)^N (\frac {1}{2}
\partial_{z}^{-1} e B)^N
\partial_{ z}^{-1}\left (e^{\tau(\Phi-\overline\Phi)}e_1
(q+\frac{\widetilde q}{\tau})\right )\vert_{G_{\epsilon_0}} =
e^{-2i\tau\psi} (\frac{m_{+,N} e^{2i\tau\psi(\widetilde x)}}
{\tau^2} + o_{C^2(\overline G_{\epsilon_0})}(\frac{1}{\tau^2})).
\end{equation}

This formula follows immediately from the stationary phase argument,
the assumption that functions $q,\widetilde q$ equal zero on
$\mathcal H$, Proposition \ref{gandon} and the representation of the
operator $(-1)^N (\frac {1}{2}
\partial_{z}^{-1} e B)^N
\partial_{ z}^{-1} e_1$ in the form:
$$
(-1)^N (\frac {1}{2}
\partial_{z}^{-1} e B)^N
\partial_{ z}^{-1}e_1 g=\int_\Omega \tilde K(x,\xi)e_1(\xi)g(\xi)d\xi,
$$
where
$$
\tilde K(x,\xi)=\frac{\tilde K_*(x,\xi)}{x_1-ix_2-(\xi_1-i\xi_2)},
\quad \tilde K_*(x,\xi)\in C^5(\bar\Omega)\times  C^5(\bar\Omega).
$$

Next let $x^0=(x^0_1,x^0_2)$ be an arbitrary fixed point in
$\Omega,$
$\partial^\beta=\partial^{\beta_1}_{x^0_1}\partial^{\beta_1}_{x^0_1},$
and $z^0=x_1^0+ix^0_2.$ Let $ \mbox{C}V=-\frac 12\partial_{ z}^{-1}
e VB$ for any matrix valued function $V(x).$
 By Proposition \ref{Proposition 3.0} there exists $\hat N$ such that the operator
\begin{equation}\label{koko}
C^N\in\mathcal L(L^\frac 43(\Omega),C^5(\overline\Omega))\quad
\forall N\ge \hat N.
\end{equation}
We write the operator $\frac{(-1)^N}{2}(\frac 12\partial_{ z}^{-1} e
B)^{ N}\partial_{ z}^{-1}$ in the form of the integral operator
$$
\frac{(-1)^N}{2}(\frac 12\partial_{ z}^{-1} e B)^{ N}\partial_{
z}^{-1}e_1g=\frac{1}{\pi}\int_\Omega\frac{{\mathcal
K}_N(x,\xi)e_1(\xi)g(\xi_1,\xi_2)}{x_1-ix_2-(\xi_1-i\xi_2)}d\xi_1d
\xi_2 .
$$
Let us estimate the kernel $\mathcal K_N.$ Observe that
\begin{equation}{\mathcal K}_N(x^0_1,x^0_2,\xi)=(-1)^N\mbox{C}^{ N}\frac{E}{2\pi (\bar z^0-\bar z)}.
\end{equation}
Since $\mbox{sup}_{z^0\in G_\epsilon,\vert\beta\vert\le 5}\Vert
\partial^\beta\frac{e}{\bar z^0-\bar z}\Vert_{L^\frac 43(\Omega)}+\mbox{sup}_{z^0\in \Omega}\Vert
\frac{1}{\bar z^0- \bar z}\Vert_{L^\frac 43(\Omega)}<\infty$ there
exists $r\in(0,1)$  independent of $N$ such that
\begin{equation}\label{k}
\mbox{sup}_{z^0\in\Omega}\Vert \partial_{ z}^{-1}\mbox{C}^{
N-1}\frac{E}{2\pi (\bar z^0-\bar z)} \Vert_{L^\frac 43(\Omega)}\le
r^{N-\hat N}.
\end{equation}
By (\ref{k}), (\ref{koko}) we obtain
\begin{equation}\label{popo1}
\Vert \mathcal K_N(x,\cdot)\Vert_{(C^5(\overline G_\epsilon)\cap
L^\infty(\Omega))\times C^5(\overline \Omega)}\le Cr^{N-\hat N}.
\end{equation}
By (\ref{popo1}) there exist a function $\mathcal K(x,\xi)\in
(C^5(\overline G_\epsilon)\cap L^\infty(\Omega))\times C^5(\overline
\Omega)$ such that
$$
\sum_{j=\hat N+2}^\infty \frac {e^{-2\tau i\psi}}{2}  (-1)^j(\frac
{1}{2}
\partial_{z}^{-1} e B)^j
\partial_{ z}^{-1}e_1 g=e^{-2\tau
i\psi}\frac 1\pi\int_\Omega\frac{\mathcal K(x,\xi)
e_1(\xi)g(\xi)}{x_1-ix_2-(\xi_1-i\xi_2)}d\xi .
$$

So, by the stationary phase argument there exists a function $m \in
(C^2(\overline G_\epsilon)\cap L^\infty(\Omega))\times C^5(\overline
\Omega)$ such that

\begin{equation}\label{momo}
\sum_{j=\hat N+2}^\infty e^{-2\tau i \psi}(-1)^j(\frac 12\partial_{
z}^{-1} e B)^{ j}\partial_{
z}^{-1}(e^{2i\tau\psi}e_1(q+\frac{\widetilde q}{\tau}))=e^{-2\tau i
\psi}( \frac{m(x)e^{2i\tau\psi(\widetilde
x)}}{\tau^2}+o_{L^\infty(G_{\widetilde
\epsilon})}(\frac{1}{\tau^2}))\,\,\forall \widetilde\epsilon >0,
\end{equation}
\begin{equation}\label{momoZ}
\sum_{j=\hat N+2}^\infty e^{-2\tau i \psi}(-1)^j(\frac 12\partial_{
z}^{-1} e B)^{ j}\partial_{
z}^{-1}(e^{2i\tau\psi}e_1(q+\frac{\widetilde q}{\tau}))=e^{-2\tau i
\psi}( \frac{m(x)e^{2i\tau\psi(\widetilde
x)}}{\tau^2}+o_{C^2(\overline G_\epsilon)}(\frac{1}{\tau^2})).
\end{equation}

By (\ref{momo}), (\ref{momoZ}), (\ref{ika}) for any positive
$\widetilde\epsilon$ we have :\begin{eqnarray}\label{A50} \frak
T_{B_1,\tau}(e_1(q+\frac{\widetilde
q}{\tau}))\vert_{G_{\widetilde\epsilon}} =
e^{\tau(\overline\Phi-\Phi)}\left (\frac{m_+ e^{2i\tau\psi
(\widetilde
x)}}{\tau^2}+o_{L^\infty(\overline G_{\widetilde\epsilon})}(\frac{1}{\tau^2})\right)\quad\mbox{as}\,\vert\tau\vert\rightarrow +\infty ,\\
\quad \frak P_{B_1,\tau}(e_1(q + \frac{\widetilde
q}{\tau}))\vert_{G_{\widetilde\epsilon}}
=e^{\tau(\Phi-\overline\Phi)}\left (\frac{m_-
e^{-2i\tau\psi(\widetilde x)}}{\tau^2} +o_{L^\infty(\overline
G_{\widetilde\epsilon})}(\frac{1}{\tau^2})\right
)\quad\mbox{as}\,\vert\tau\vert\rightarrow +\infty.
\end{eqnarray}
and

\begin{eqnarray}\label{AA50} \frak
T_{B_1,\tau}(e_1(q+\frac{\widetilde q}{\tau}))\vert_{G_{\epsilon}} =
e^{\tau(\overline\Phi-\Phi)}\left (\frac{m_+ e^{2i\tau\psi
(\widetilde
x)}}{\tau^2}+o_{C^2(\overline G_{\epsilon})}(\frac{1}{\tau^2})\right)\quad\mbox{as}\,\vert\tau\vert\rightarrow +\infty ,\\
\quad \frak P_{B_1,\tau}(e_1(q + \frac{\widetilde
q}{\tau}))\vert_{G_{\epsilon}} =e^{\tau(\Phi-\overline\Phi)}\left
(\frac{m_- e^{-2i\tau\psi(\widetilde x)}}{\tau^2} +o_{C^2(\overline
G_{\epsilon})}(\frac{1}{\tau^2})\right
)\quad\mbox{as}\,\vert\tau\vert\rightarrow +\infty.
\end{eqnarray}

Let positive $\hat \epsilon $ be such that $\mbox{supp}\,
(1-e)\subset G_{\hat \epsilon}$ and $\hat\epsilon<\epsilon,
\epsilon''\in(\hat \epsilon,\epsilon)$ Then using (\ref{A50}) we
have
\begin{eqnarray}\label{solonka}
e^{-2\tau i\psi} T_{B_1}(e^{\tau (\Phi-\overline\Phi)}(1-e)B_1 \frak
T_{B_1,\tau}e_1(q+\frac{\widetilde q}{\tau}) )= e^{-2\tau i\psi} T_{
B_1}((1-e)B_1\frak T_{B_1} (e^{\tau
(\Phi-\overline\Phi)}e_1(q+\frac{\widetilde q}{\tau}))) \nonumber\\
= e^{-2\tau i\psi+2i\tau\psi(\widetilde x)} T_{
B_1}((1-e)\chi_{G_{\epsilon ''}}B_1 \frac{m_+
 }{\tau^2} +(1-e)\chi_{G_{\epsilon ''}}o_{C^2( \overline{G_{\epsilon''}})}(\frac{1}{\tau^2}))+\nonumber\\
 e^{-2\tau i\psi+2i\tau\psi(\widetilde x)} T_{ B_1}((1-e)(1-\chi_{G_{\epsilon ''}})B_1
\frac{m_+
 }{\tau^2} +(1-e)(1-\chi_{G_{\epsilon ''}})o_{L^\infty( G_{\hat \epsilon})}(\frac{1}{\tau^2})).
\end{eqnarray}  Here in order to obtain the last equality we
used (\ref{A50}) and (\ref{short}). Using (\ref{AA50}),
(\ref{solonka}), (\ref{sonka}) and Proposition \ref{Proposition 3.0}
we obtain (\ref{50}). $\blacksquare$

Denote $ q_1=P_{A_1}((Q_1-2\partial_z A_1-B_1A_1)w_{0})-M_1, q_2=
T_{B_1}((Q_1-2\partial_{\overline z} B_1-A_1B_1)\widetilde
w_{0})-M_2\in C^{5+\alpha}(\bar \Omega)$, where the functions
$M_1\in Ker (\partial_{\overline z}+A_1)$ and $M_2\in Ker
(\partial_z+B_1) $ are taken such that
\begin{equation}\label{kl}
q_1(\widetilde x)=q_2(\widetilde x)=0, \quad
\partial^{\alpha_1}_{x_1}\partial^{\alpha_2}_{x_2} q_1
(x)=\partial^{\alpha_1}_{x_1}\partial^{\alpha_2}_{x_2} q_2(x)\quad
\forall x\in \mathcal H\setminus\{\widetilde x\}\,\,\,
\mbox{and}\,\,\, \forall \alpha_1+\alpha_2\le 5.
\end{equation}
By Proposition \ref{popo}, there exist functions $m_\pm\in
C^2(\partial\Omega)$ such that
\begin{equation}\label{50l}
\widetilde{\mathcal R}_{\tau, B_1}(e_1(q_1+\frac{\widetilde q_1}{\tau})) =
e^{\tau(\overline\Phi-\Phi)}\left (\frac{m_+ e^{2i\tau\psi (\widetilde
x)}}{\tau^2}+o_{H^1(\partial\Omega)}(\frac{1}{\tau^2})\right)
\quad\mbox{as}\,\vert\tau\vert\rightarrow +\infty
\end{equation}
and
\begin{equation}\label{50ll}
\quad{\mathcal R}_{\tau, A_1} (e_1(q_2 + \frac{\widetilde q_2}{\tau}))
=e^{\tau(\Phi-\overline\Phi)}
\left (\frac{m_- e^{-2i\tau\psi(\widetilde
x)}}{\tau^2} +o_{H^1(\partial\Omega)}(\frac{1}{\tau^2})\right
)\quad\mbox{as}\,\vert\tau\vert\rightarrow +\infty.
\end{equation}

Next we introduce the functions $w_{-1},\widetilde w_{-1},
a_\pm,b_\pm\in C^2(\overline \Omega)$ as a solutions to the following
boundary value problems:
\begin{equation}\label{zad-1}
\mathcal K(x,D)(w_{-1},\widetilde w_{-1})=0\quad\mbox{in}\,\,\Omega,\quad
(w_{-1}+\widetilde
w_{-1})\vert_{\Gamma_0}=\frac{q_1}{2\partial_z\Phi}
+ \frac{q_2}{2\partial_{\overline z}\overline\Phi},
\end{equation}
$$
\partial^{\alpha_1}_{x_1}\partial^{\alpha_2}_{x_2} w_{-1}(x)
=\partial^{\alpha_1}_{x_1}\partial^{\alpha_2}_{x_2} \widetilde
w_{-1}(x)\quad \forall x\in \mathcal H\,\,\, \mbox{and}\,\,\,
\forall \alpha_1+\alpha_2\le 2,
$$
\begin{equation}
\mathcal K(x,D)(a_\pm,b_\pm)=0\quad\mbox{in}\,\,\Omega,\quad
(a_\pm+b_\pm)\vert_{\Gamma_0}=m_\pm .
\end{equation}

We set $p_1=-(Q_1-2\partial_{\overline
z}B_1-A_1B_1)(\frac{e_1q_1}{2\partial_z\Phi}+w_{-1})+L_1(x,D)
(\frac{e_2q_1}{2\partial_z\Phi})$, $p_2=-(Q_1-2\partial_{
z}A_1-B_1A_1)(\frac{e_1q_2}{2\partial_{\overline z}\overline\Phi}+\widetilde
w_{-1})+L_1(x,D) (\frac{e_2q_2}{2\partial_{\overline z}\overline\Phi}),$
$\widetilde q_2=T_{B_1}p_2-\widetilde M_2, \widetilde q_1=P_{A_1}p_1
-\widetilde M_1$, where  $\widetilde M_1\in Ker (\partial_{\overline z}+A_1)$
and $\widetilde M_2\in Ker (\partial_z+B_1)$ are taken such that
\begin{equation}\label{gandon1}
\widetilde q_1(\widetilde x)
=\widetilde q_2(\widetilde x)=0,\quad \partial^{\alpha_1}_{x_1}
\partial^{\alpha_2}_{x_2} \widetilde q_1
(x)=\partial^{\alpha_1}_{x_1}\partial^{\alpha_2}_{x_2} \widetilde
q_2(x)\quad \forall x\in \mathcal H\setminus\{\widetilde x\}\,\,\,
\mbox{and}\,\,\, \forall \alpha_1+\alpha_2\le 2.
\end{equation}
Since $\frac{\widetilde q_1}{2\partial_z\Phi},\frac{\widetilde q_2}
{2\partial_{\overline z}\overline\Phi}\in H^1(\partial\Omega)$
by (\ref{gandon1}), there exists a solution $(w_{-2},\widetilde w_{-2})
\in H^1(\overline\Omega)$ to the
boundary value problem
\begin{equation}\label{zad-1} \mathcal K(x,D)(w_{-2},\widetilde w_{-2})
=0\quad\mbox{in}\,\,\Omega,\quad
(w_{-2}+\widetilde w_{-2})\vert_{\Gamma_0}=\frac{\widetilde
q_1}{2\partial_z\Phi}+\frac{\widetilde q_2}{2\partial_{\overline z}
\overline\Phi}.
\end{equation}

We introduce the functions $w_{0,\tau}, \widetilde w_{0,\tau}
\in H^1(\Omega)$ by
\begin{equation}\label{zad1}
w_{0,\tau}=w_0+\frac{w_{-1}-e_2q_1/2\partial_{
z}\Phi}{\tau}+\frac{1}{\tau^2}( e^{2i\tau\psi(\widetilde
x)}a_++e^{-2i\tau\psi(\widetilde x)}a_-+w_{-2}-\frac{\widetilde
q_1}{2\partial_z\Phi})
\end{equation}
and
\begin{equation}\label{zad2}
\widetilde w_{0,\tau}=\widetilde w_0+\frac{\widetilde w_{-1}-e_2
q_2/2\partial_{\overline z}\overline\Phi}{\tau}+\frac{1}{\tau^2}(
e^{2i\tau\psi(\widetilde x)}b_++e^{-2i\tau\psi(\widetilde
x)}b_-+\widetilde w_{-2}-\frac{\widetilde q_2}{2\partial_{\overline z}
\overline\Phi}).
\end{equation}

Simple computations and Proposition 8  for any $p\in (1,\infty)$
imply the asymptotic formula:
\begin{eqnarray}\label{251}
L_1(x,D)(-e^{\tau\Phi}\widetilde{\mathcal R}_{\tau, B_1}(e_1(q_1+\widetilde
q_1/\tau))-\frac{e_2(q_1+\widetilde
q_1/\tau)e^{\tau\Phi}}{2\tau\partial_z\Phi}-e^{\tau\overline\Phi}{\mathcal
R}_{\tau,
A_1}(e_1(q_2+\widetilde q_2/\tau))\nonumber\\-\frac{e_2(q_2+\widetilde q_2/\tau)e^{\tau\overline\Phi}}{2\tau\partial_{\overline z}\overline\Phi})
                                                      \nonumber
= -L_1(x,D)(e^{\tau\Phi}\widetilde{\mathcal R}_{\tau, B_1}
(e_1(q_1+\widetilde q_1/\tau))+\frac{e_2(q_1+\widetilde
q_1/\tau)e^{\tau\Phi}}{2\tau\partial_z\Phi})\\
- L_1(x,D)(e^{\tau\overline\Phi}{\mathcal R}_{\tau,
A_1}(e_1(q_2+\widetilde q_2/\tau))+\frac{e_2(q_2+\widetilde q_2/\tau)
e^{\tau\overline\Phi}}{2\tau\partial_{\overline z}\overline\Phi})\nonumber\\
= -(Q_1-2\partial_{\overline z}B_1-A_1B_1)e^{\tau\Phi}\widetilde {\mathcal
R}_{\tau,B_1}(e_1(q_1+\widetilde q_1/\tau))-(Q_1-2\partial_{
z}A_1-B_1A_1)e^{\tau\overline\Phi}{\mathcal
R}_{\tau,A_1}(e_1(q_2+\widetilde q_2/\tau))\nonumber\\
-e^{\tau\Phi}L_1(x,D) (\frac{e_2(q_1+\widetilde
q_1/\tau)}{2\tau\partial_z\Phi})-e^{\tau\overline\Phi}L_1(x,D)
(\frac{e_2(q_2+\widetilde q_2/\tau)}{2\tau\partial_{\overline z}\overline\Phi})
\nonumber\\-(Q_1-2\partial_{\overline z} B_1-A_1B_1)\widetilde w_0
e^{\tau\overline\Phi}-(Q_1-2\partial_z
A_1-B_1A_1)w_0e^{\tau\Phi}\nonumber \\
-(Q_1-2\partial_{\overline z}B_1-A_1B_1)e^{\tau\Phi}\frac{e_1q_1}
{2\tau\partial_z\Phi}-(Q_1-2\partial_{z}A_1-B_1A_1)
e^{\tau\overline\Phi}\frac{e_1q_2}{2\tau\partial_{\overline z}
\overline\Phi}\nonumber\\
+ \frac 1\tau((Q_1-2\partial_{\overline
z}B_1-A_1B_1)\frac{e_1q_1}{2\partial_z\Phi}+L_1(x,D)
(\frac{e_2q_1}{2\partial_z\Phi}))e^{\tau\Phi} \nonumber\\
+ \frac 1\tau((Q_1-2\partial_{
z}A_1-B_1A_1)\frac{e_2q_2}{2\partial_{\overline
z}\overline\Phi}+L_1(x,D) (\frac{e_2q_2}{2\partial_{\overline
z}\overline\Phi}))e^{\tau\overline\Phi}\nonumber\\
= - \frac{1}{\tau}(Q_1-2\partial_{\overline
z}B_1-A_1B_1) w_{-1}e^{\tau\Phi}-\frac{1}{\tau}(Q_1-2\partial_{z}
A_1-B_1A_1)\widetilde
w_{-1}e^{\tau\overline\Phi}\nonumber\\
- (Q_1-2\partial_{\overline z}
B_1-A_1B_1)\widetilde w_0 e^{\tau\overline\Phi}-(Q_1-2\partial_z
A_1-B_1A_1)w_0e^{\tau\Phi}
+e^{\tau\varphi}o_{L^p(\Omega)}(\frac{1}{\tau}).
\end{eqnarray}

Using this formula, we prove the following proposition.
\begin{proposition}\label{Proposition 00}
For any $p>1$, we have the asymptotic formula:
\begin{eqnarray}\label{249}
L_1(x,D)(w_{0,\tau}e^{\tau \Phi}+\widetilde w_{0,\tau}
e^{\tau \overline
\Phi}-e^{\tau\Phi}\widetilde{\mathcal R}_{\tau, B_1}
(e_1(q_1+\widetilde
q_1/\tau))-e^{\tau\overline\Phi}{\mathcal R}_{\tau,
A_1}(e_1(q_2+\widetilde
q_2/\tau)))=e^{\tau\varphi}o_{L^p(\Omega)}(\frac{1}{\tau}) ,\\
\quad (w_{0,\tau}e^{\tau \Phi}+\widetilde w_{0,\tau} e^{\tau \overline
\Phi}-e^{\tau\Phi}\widetilde{\mathcal R}_{\tau, B_1}(e_1(q_1+\widetilde
q_1/\tau))-e^{\tau\overline\Phi}{\mathcal R}_{\tau,
A_1}(e_1(q_2+\widetilde
q_2/\tau)))\vert_{\Gamma_0}=e^{\tau\varphi}o_{H^1(\Gamma_0)}
(\frac {1}{\tau^2}).
\end{eqnarray}
\end{proposition}

{\bf Proof.} By (\ref{zzz}), (\ref{50l}), (\ref{50ll}),
(\ref{zad-1}) and (\ref{zad-1})-(\ref{zad2}), we have
$$
(w_{0,\tau}e^{\tau \Phi}+\widetilde w_{0,\tau} e^{\tau \overline
\Phi}-e^{\tau\Phi}\widetilde{\mathcal R}_{\tau, B_1}(e_1(q_1+\widetilde
q_1/\tau))-e^{\tau\overline\Phi}{\mathcal R}_{\tau,
A_1}(e_1(q_2+\widetilde q_2/\tau))\vert_{\Gamma_0}
$$
$$
= (w_{0,\tau}e^{\tau \varphi}+\widetilde w_{0,\tau} e^{\tau
\varphi}-e^{\tau\varphi}\widetilde{\mathcal R}_{\tau,
B_1}(e_1(q_1+\widetilde q_1/\tau))-e^{\tau\varphi}{\mathcal R}_{\tau,
A_1}(e_1(q_2+\widetilde q_2/\tau)))\vert_{\Gamma_0}
$$
$$
= e^{\tau\varphi}(w_0+\frac{w_{-1}-e_2q_1/2\partial_{ z}\Phi}{\tau}
+\frac{1}{\tau^2}(
e^{2i\tau\psi(\widetilde x)}a_++e^{-2i\tau\psi(\widetilde x)}a_-
+w_{-2}-\frac{\widetilde q_1}{2\partial_z\Phi}) $$
$$
+ \widetilde w_0+\frac{\widetilde w_{-1}-e_2 q_2/2\partial_{\overline
z}\overline\Phi}{\tau}+\frac{1}{\tau^2}( e^{2i\tau\psi(\widetilde
x)}b_+
+e^{-2i\tau\psi(\widetilde x)}b_-+\widetilde w_{-2}
-\frac{\widetilde q_2}{2\partial_{\overline z}\overline\Phi})
$$
$$
-e^{\tau\varphi}\widetilde{\mathcal R}_{\tau, B_1}(e_1(q_1+\widetilde
q_1/\tau))-e^{\tau\varphi}{\mathcal R}_{\tau, A_1}(e_1(q_2+\widetilde
q_2/\tau))))\vert_{\Gamma_0}
$$
$$
= e^{\tau\varphi}\{\frac{1}{\tau^2}
(e^{2i\tau\psi(\widetilde x)}a_++e^{-2i\tau\psi(\widetilde
x)}a_- + e^{2i\tau\psi(\widetilde
x)}b_++e^{-2i\tau\psi(\widetilde x)}b_-)
$$
$$
- e^{\tau\varphi}\widetilde{\mathcal R}_{\tau, B_1}(e_1(q_1+\widetilde
q_1/\tau))-e^{\tau\varphi}{\mathcal R}_{\tau, A_1}(e_1(q_2+\widetilde
q_2/\tau))\}\vert_{\Gamma_0}=e^{\tau\varphi}o_{H^1(\Gamma_0)}(\frac
{1}{\tau^2}).
$$
Here in order to obtain the final equality, we used Proposition
\ref{popo}. Similarly to (\ref{oi}) we obtain
\begin{eqnarray}\label{250}
L_1(x,D)(w_{0,\tau}e^{\tau\Phi}+\widetilde
w_{0,\tau}e^{\tau\overline\Phi}-\frac{e_2(q_1+\widetilde
q_1/\tau)e^{\tau\Phi}}{2\tau\partial_z\Phi}-\frac{e_2(q_2+\widetilde
q_2/\tau)e^{\tau\overline\Phi}}{2\tau\partial_{\overline
z}\overline\Phi})\nonumber\\
= (Q_1-2\partial_{\overline z}B_1-A_1B_1) (w_{0,\tau}-\frac{e_2(q_1+\widetilde
q_1/\tau)}{2\tau\partial_z\Phi})e^{\tau\Phi} \nonumber\\
+ (Q_1-2\partial_{z}A_1-B_1A_1)(\widetilde w_{0,\tau}
-\frac{e_2(q_2+\widetilde
q_2/\tau)}{2\tau\partial_{\overline z}\overline\Phi})e^{\tau\overline\Phi}.
\end{eqnarray}
By (\ref{250}) and (\ref{251}), we obtain (\ref{249}). $\blacksquare$

We set ${\mathcal O}_{\epsilon}=\{x\in \Omega; \thinspace dist (x,\partial
\Omega)\le \epsilon\}.$
In order to construct the last term in complex geometric optics solution,
we need the following proposition:

\begin{proposition}\label{Proposition 0}
Let $A,B\in C^{5+\alpha}(\overline\Omega)$ and $Q\in C^{4+\alpha}
(\overline\Omega)$ for some $\alpha\in (0,1),$ $f\in L^p(\Omega)$
for some $p>2$, $dist (\overline\Gamma_0,supp\,f)>0$, $q\in H^\frac
12(\Gamma_0),$ and $\epsilon$ be a small positive number such that $
\overline{{\mathcal O}_{\epsilon}}\cap({\mathcal
H}\setminus\Gamma_0)=\emptyset.$ Then there exists $C$ independent
of $\tau$ and $\tau_0$ such that for all $\vert\tau\vert>\tau_0$,
there exists a solution to the boundary value problem
\begin{equation}\label{mimino}
L(x,D)w=fe^{\tau\Phi}\quad \mbox{in}\,\,\Omega, \quad
w\vert_{\Gamma_0} =qe^{\tau\varphi}/\tau
\end{equation}
such that
\begin{equation}\label{mimino111}
\root\of{\vert\tau\vert} \Vert w e^{-\tau\varphi}\Vert_{L^2(\Omega)}
+ \frac{1}{\root\of{\vert\tau\vert} } \Vert (\nabla w)
e^{-\tau\varphi}\Vert_{L^2(\Omega)}+\Vert we^{-\tau\varphi}
\Vert_{H^{1,\tau}(\mathcal O_{\epsilon})} \le C(\Vert
f\Vert_{L^p(\Omega)} + \Vert q \Vert_{H^\frac 12(\Gamma_0)}).
\end{equation}
\end{proposition}

{\bf Proof.} First let us assume that $f$ is identically equal to
zero. Let $(d,\widetilde d)\in H^1(\Omega)\times  H^1(\overline\Omega)$
satisfy
\begin{equation}\label{zadnitsa}
\mathcal K(x,D)(d,\widetilde d)=0\quad\mbox{in}\,\Omega, \quad (d+\widetilde
d)\vert_{\Gamma_0}=q.
\end{equation}
For exitance of such a solution see e.g. \cite{Wendland}. By
(\ref{oi}) and (\ref{zadnitsa}), we have
$$
L(x,D)(\frac{d}{\tau}e^{\tau\Phi}+\frac{\widetilde
d}{\tau}e^{\tau\overline\Phi})=\frac{1}{\tau}(Q-2\partial_z A-BA)d
e^{\tau\Phi}+\frac{1}{\tau}(Q-2\partial_{\overline z} B-AB)\widetilde d
e^{\tau\overline\Phi}.
$$
By Proposition \ref{vanka}, there exists  a solution $w$ to the
boundary value problem
$$
L(x,D)\tilde w=-\frac{1}{\tau}(Q-2\partial_z A-BA)d
e^{\tau\Phi}-\frac{1}{\tau}(Q-2\partial_{\overline z}
B-AB)\widetilde d e^{\tau\overline\Phi},\quad \tilde
w\vert_{\Gamma_0}=0
$$
such that there exists a constant $C>0$ such that
$$
\Vert \tilde we^{-\tau\varphi}\Vert_{H^{1,\tau}(\Omega)}\le \frac{C}
{\root\of{\vert \tau\vert}}\Vert (Q-2\partial_z A-BA)d
e^{i\tau\psi}+(Q-2\partial_{\overline z} B-AB)\widetilde d
e^{-\tau\psi})\Vert_{L^2(\Omega)}\le \frac{C
}{\root\of{\vert\tau\vert}}\Vert q\Vert_{H^\frac 12(\Gamma_0)}
$$
for all large $\tau>0$.

Then the function $(\frac{d}{\tau}e^{\tau\Phi}+\frac{\widetilde
d}{\tau}e^{\tau\overline\Phi})+\tilde w$ is a solution to
(\ref{mimino}) which satisfies (\ref{mimino111}) if $f\equiv 0.$

If $f$ is not identically equal zero, then we consider the function
$\tilde w=\widetilde e e^{\tau\Phi}\widetilde{\mathcal R}_{\tau,
B}(e_1q_0)$, where $\widetilde e\in
C^\infty_0(\Omega),\,\,\widetilde e\vert_{\mbox{supp}e_1}=1$ and
$q_0=P_{A}f-M,$  where a function $M\in C^5(\bar \Omega)$ belongs to
$Ker\, (2\partial_z+B)$ and chosen such that $q_0\vert_{\mathcal
H}=0.$ Then $ L(x,D)\tilde w=(Q-2\partial_{\overline z}B-AB)\tilde
w+ \widetilde ee_1fe^{\tau\Phi}+2\widetilde
ee^{\tau\Phi}q_0\partial_{\overline
z}e_1+e^{\tau\Phi}(2\partial_{\overline z}+A)(\partial_z
e\widetilde{\mathcal R}_{\tau, B}(e_1q_0)). $ Since, by Proposition
8, the function $\widetilde f(\tau,\cdot)=e^{-\tau\Phi}L(x,D)\tilde
w-f$ can be represented as a sum of two functions, where the first
one equal to zero in a neighborhood of $\mathcal H$ and is bounded
uniformly in $\tau$ in $L^2(\Omega)$ norm, the second one is
$O_{L^2(\Omega)}(\frac{1}{\tau})$. Applying Proposition \ref{vanka}
to the boundary value problem
$$
L(x,D)w_*=\widetilde fe^{\tau\Phi}\quad\mbox{in}\,\,\Omega, \quad
w_*\vert_{\Gamma_0}=0,
$$
we construct a solution such that
$$
\Vert w_{*}e^{-\tau\varphi}\Vert_{H^{1,\tau}(\Omega)}\le C\Vert
f\Vert_{L^p(\Omega)}.
$$ The function $w^*-\tilde w$ solves the boundary value problem
(\ref{mimino}) and satisfies estimate $(\ref{mimino111}).$
 $\blacksquare$

Using Propositions \ref{Proposition 0} and \ref{Proposition 00},
we construct the last term $u_{-1}$ in complex geometric
optics solution which satisfies
\begin{equation}\label{mimino11}
\root\of{\vert\tau\vert} \Vert u_{-1} \Vert_{L^2(\Omega)} +
\frac{1}{\root\of{\vert\tau\vert} } \Vert (\nabla u_{-1})
\Vert_{L^2(\Omega)}+\Vert u_{-1} \Vert_{H^{1,\tau}(\mathcal
O_{\epsilon})}=o(\frac{1}{\tau})\quad \mbox{as}\,\,\tau\rightarrow
+\infty.
\end{equation}
Finally we obtain a complex geometric optics solution in the form:
\begin{equation}\label{zad}
u_1(x)=w_{0,\tau}e^{\tau \Phi}+\widetilde w_{0,\tau} e^{\tau \overline
\Phi}-e^{\tau\Phi}\widetilde{\mathcal R}_{\tau, B_1}(q_1+\widetilde
q_1/\tau)-e^{\tau\overline\Phi}{\mathcal R}_{\tau, A_1}(q_2+\widetilde
q_2/\tau)+e^{\tau \varphi} u_{-1}.
\end{equation}

Obviously
\begin{equation}\label{zad22}
L_1(x,D)u_1=0\quad\mbox{in}\,\,\Omega,\quad u_1\vert_{\Gamma_0}=0.
\end{equation}

Let $u_1$ be a complex geometrical optics  solution as in
(\ref{zad}).

Let $\mbox{\bf e}\in C^\infty_0(\Bbb R^n)$ be a function such that
$\mbox{\bf e}$ is equal to one in a ball of small radius centered at
$0.$ We set
\begin{equation}\label{BIGzopa}
\eta(x,s)=\mbox{\bf e}((x-\widetilde x)e^{s^2}).
\end{equation}

Then the operator
\begin{eqnarray*}
&&L_2(x,s,D)=e^{-s\eta}L_2(x,D)e^{s\eta}=\Delta
+2(A_2+2s\eta_{\overline z})\partial_z+2(B_2+2s\eta_{ z})
\partial_{\overline z}+Q_2\\
+ &&(s\Delta\eta +s^2(\nabla\eta,\nabla\eta))E+2s\eta_z
A_2+2s\eta_{\overline z} B_2
\end{eqnarray*}
is of the form (\ref{o-1}) and has the
same partial Cauchy data as the operator $L_2(x,D)$.
Also for the operator $L_2(x,s,D)$, one can construct a similar complex
geometric optics  solution.

Consider the operator
\begin{eqnarray}
{L}_2(x,s,{D})^{*} =4\partial_{ z} {\partial}_{\overline z}
-2{{A}^*_{2,s}}{\partial}_{\overline z}
-2{{B}^*_{2,s}}\partial_z+{Q^*_{2,s}} -2\partial_{\overline z}
{A}_{2,s}^*
-2\partial_z {B}^*_{2,s}\nonumber\\
=(2\partial_z-{{A}^*_{2,s}}) (2\partial_{\overline z}-{{B}^*_{2,s}})
+Q_{2}^* -2
\partial_{\overline z} {{A}^*_{2}}
-{{A}^*_{2}}  {{B}^*_{2}}\nonumber\\
=(2\partial_{\overline z}-{{B}^*_{2,s}}) (2\partial_z-{{A}^*_{2,s}})
+Q^*_{2} -2
\partial_z {{B}^*_{2}} -{{B}^*_{2}} {{A}^*_{2}}.
                              \nonumber
\end{eqnarray}
Similarly we construct the complex geometric optics solutions to
the operator $L_2(x,s,D)^*.$
Let $(w_1,\widetilde w_1) \in C^{6+\alpha}(\overline \Omega)$
be a solutions to the following boundary value problem:
\begin{equation}\label{ll1}
\mathcal M(x,D)(w_{1},\widetilde w_{1})=((2\partial_{\overline
z}-{B_2^*})w_{1},(2\partial_{ z}-{A_2^*})\widetilde
w_{1})=0\quad\mbox{in}\,\,\Omega, \quad  (w_{1}+\widetilde
w_{1})\vert_{\Gamma_0}=0,
\end{equation}
$$
\partial^{\alpha_1}_{x_1}\partial^{\alpha_2}_{x_2} w_1(x)
=\partial^{\alpha_1}_{x_1}\partial^{\alpha_2}_{x_2} \widetilde
w_1(x)\quad \forall x\in \mathcal H\,\,\, \mbox{and}\,\,\, \forall
\alpha_1+\alpha_2\le 2, $$
\begin{equation}\label{iiii}
\lim_{x\rightarrow x_\pm}\frac{ \vert w_1(x)\vert}{\vert x-x_\pm\vert^{98}}
=\lim_{x\rightarrow x_\pm}\frac{ \vert \widetilde w_1(x)\vert}
{\vert x-x_\pm\vert^{98}}=0.
\end{equation}

Such a pair $(w_1,\widetilde w_1)$ exists due
to Proposition 9.
 We set $(w_{1,s},\widetilde
w_{1,s})=e^{s\eta}(w_{1},\widetilde w_{1}).$ Observe that
$$
{L}_2(x,s,{D})^{*}(w_{1,s}e^{-\tau\Phi})=(Q_{2}^* -2
\partial_{\overline z} {{A}^*_{2}}-{{A}^*_{2}} {{B}^*_{2}}
)w_{1,s}e^{-\tau\Phi},
$$
$$
{L}_2(x,s,{D})^{*}(\widetilde
w_{1,s}e^{-\tau\overline\Phi})=(Q_{2}^* -2 \partial_z {{B}^*_{2}}
-{{B}^*_{2}} {{A}^*_{2}})\widetilde w_{1,s} e^{-\tau\overline\Phi}.
$$
We set
\begin{equation}\label{loko}
P_{- B^*_{2,s}}=e^{s\eta}P_{- B^*_{2}}e^{-s\eta}, T_{-
A^*_{2,s}}=e^{s\eta}T_{- A^*_{2}}e^{-s\eta},\widetilde{\mathcal
R}_{-\tau,- A_{2,s}^*}=e^{s\eta}\widetilde{\mathcal R}_{-\tau,-
A_{2}^*}e^{-s\eta},\quad {\mathcal R}_{-\tau,-
B_{2,s}^*}=e^{s\eta}{\mathcal R}_{-\tau,- B_{2}^*}e^{-s\eta},
\end{equation}
\begin{equation}\label{loko1}
q_3=P_{- B^*_{2}}(( Q_{2}^* -2
\partial_{\overline z} {{A}_2}^* -{{A}_{2}}^* {{B}_{2}}^*
)w_{1})-M_{3},\quad q_4=T_{- A^*_{2}}(( Q_2^* -2
\partial_z {{B}^*_2} -{{B}_{2}}^* {{A}_2}^*)\widetilde w_{1})-M_{4}.
\end{equation}
Denote $q_{3,s}=P_{-B^*_{2,s}}(( Q_{2}^* -2
\partial_{\overline z} {{A}_2}^* -{{A}_{2}}^* {{B}_{2}}^* )w_{1,s})-M_{3,s}
=e^{s\eta} q_3, q_{4,s}=T_{- A^*_{2,s}}(( Q_2^* -2
\partial_z {{B}^*_2} -{{B}_{2}}^* {{A}_2}^*)\widetilde w_{1,s})-M_{4,s}
=e^{s\eta} q_4 $ where the functions $M_{j,s}=e^{s\eta}M_j,$ $M_{3}\in Ker
(2\partial_{\overline z}-B_{2}^*)$ and
$M_{4}\in Ker (2\partial_{ z}-A_{2}^*)$ are chosen such that
\begin{equation}\label{lada1}
q_3(\widetilde x)=q_4(\widetilde x)=0, \quad
\partial^{\alpha_1}_{x_1}\partial^{\alpha_2}_{x_2} q_3(x)
=\partial^{\alpha_1}_{x_1}\partial^{\alpha_2}_{x_2} q_4(x)\quad
\forall x\in \mathcal H\setminus\{\widetilde x\}\,\,\,
\mbox{and}\,\,\, \forall \alpha_1+\alpha_2\le 5.
\end{equation}

By (\ref{lada1}) the functions
$\frac{q_3}{2\partial_z\Phi},\frac{q_4}{2\partial_{\overline
z}\overline\Phi}$ belong to the space $C^2(\overline\Gamma_0).$
 Therefore  we  can
introduce the functions $w_{-3}, \widetilde w_{-3}, \widetilde
 a_\pm, \widetilde  b_\pm \in C^2(\overline \Omega)$ as a solutions
to the following boundary value problems:
\begin{equation}
\mathcal M(x,D)(w_{-3},\widetilde w_{-3})=0\quad\mbox{in}\,\,\Omega,\quad
(w_{-3}+\widetilde
w_{-3})\vert_{\Gamma_0}=\frac{q_3}{2\partial_z\Phi}+\frac{q_4}
{2\partial_{\overline z}\overline\Phi},
\end{equation}
$$
\partial^{\alpha_1}_{x_1}\partial^{\alpha_2}_{x_2} w_{-3}(x)
= \partial^{\alpha_1}_{x_1}\partial^{\alpha_2}_{x_2} \widetilde
w_{-3}(x)\quad \forall x\in \mathcal H \,\,\, \mbox{and}\,\,\,
\forall \alpha_1+\alpha_2\le 2,
$$
\begin{equation}
\mathcal M(x,D)(\widetilde a_\pm,\widetilde
b_\pm)=0\quad\mbox{in}\,\,\Omega,\quad (\widetilde a_\pm+\widetilde
b_\pm)\vert_{\Gamma_0}=\widetilde m_\pm .
\end{equation}
Let
$$
p_3=-(Q_{2}^*-2\partial_{\overline
z}A_{2}^*-A_{2}^*B_{2}^*)(\frac{e_1q_{3,s}}{2\partial_z\Phi}
+w_{-3,s})-L_2(x,s,D)^*(\frac{q_{3,s}e_2}{2\partial_z\Phi}),
$$
$$
p_4=-(Q_{2}^*-2\partial_zB_{2}^*-B_{2}^*A_{2}^*)
(\frac{e_1q_{4,s}}{2\partial_{\overline
z}\overline\Phi}+\widetilde w_{-3,s})
-L_2(x,s,D)^*(\frac{q_{4,s}e_2}{2\partial_{\overline
z}\overline\Phi})
$$
and
$$
\widetilde q_3=e^{-s\eta}(P_{-
B^*_{2,s}}p_3-\widetilde M_{3,s}), \quad
\widetilde q_4=e^{-s\eta}(T_{-A^*_{2,s}} p_4-\widetilde M_{4,s}),
$$
where $\widetilde M_{3,s}\in Ker (2\partial_{\overline z}-B_{2,s}^*),
\widetilde M_{4,s}\in Ker
(2\partial_{ z}-A_{2,s}^*),$ and $(\widetilde q_{3,s},\widetilde
q_{4,s})=e^{s\eta}(\widetilde q_3,\widetilde q_4)$ are chosen such that
 \begin{equation}\label{lada}
\widetilde q_{3,s}(\widetilde x)=\widetilde q_{4,s}(\widetilde x)=0,\quad
\partial^{\alpha_1}_{x_1}\partial^{\alpha_2}_{x_2} \widetilde q_{3,s}(x)
= \partial^{\alpha_1}_{x_1}\partial^{\alpha_2}_{x_2} \widetilde
q_{4,s}(x)\quad \forall x\in \mathcal H\setminus\{\widetilde
x\}\,\,\, \mbox{and}\,\,\, \forall \alpha_1+\alpha_2\le 2.
\end{equation}

The following asymptotic formula holds true:
\begin{proposition} \label{lp}
There exist smooth functions $\widetilde m_\pm\in C^2(\partial\Omega)$,
independent of $\tau$ and $s$, such that
\begin{equation}
\widetilde{\mathcal R}_{-\tau,- A^*_{2,s}}(e_1(q_{3,s}+\widetilde
q_{3,s}/\tau))=\frac{\widetilde m_+ e^{2i\tau(\psi+\psi(\widetilde
x))}}{\tau^2}+e^{2i\tau\psi}o_{H^1(\partial
\Omega)}(\frac{1}{\tau^2})\quad\mbox{as}\,\,\vert\tau\vert\rightarrow
+\infty
\end{equation}
and
\begin{equation}
{\mathcal R}_{-\tau,- B_{2,s}^*}(e_1(q_{4,s}+\widetilde
q_{4,s}/\tau))=\frac{\widetilde m_- e^{-2i\tau(\psi+\psi(\widetilde
x))}}{\tau^2}+e^{-2i\tau\psi}o_{H^1(\partial
\Omega)}(\frac{1}{\tau^2})\quad\mbox{as}\,\,\vert\tau\vert\rightarrow
+\infty.
\end{equation}
\end{proposition}

{\bf Proof.} The functions $q_{3,s}, q_{4,s}$ belong to the space
$C^{5+\alpha}(\bar \Omega)$ $q_{3,s}, q_{4,s}$ belong to the space
$W^1_p( \Omega)$ for any $p>1$.
By (\ref{lada1}) and (\ref{lada}), we have
$q_{3,s}=q_{4,s}=\widetilde q_{3,s}=\widetilde q_{4,s}=0$ on
$\mathcal H.$  By (\ref{loko}) and (\ref{loko1}), we have
$$
\widetilde{\mathcal R}_{-\tau,- A_{2,s}^*}(e_1(q_{3,s}+\widetilde
q_{3,s}/\tau))=e^{s\eta}\widetilde{\mathcal R}_{-\tau,-
A_{2}^*}(e_1(q_{3}+\widetilde q_{3}/\tau))
$$
and
$$
{\mathcal R}_{-\tau,-B_{2,s}^*}(e_1(q_{4,s}+\widetilde q_{4,s}/\tau))
=e^{s\eta}{\mathcal
R}_{-\tau,- B_{2}^*}(e_1(q_{4}+\widetilde q_{4}/\tau)).
$$
Then applying Proposition \ref{popo} and taking into account
(\ref{BIGzopa}), we obtain Proposition \ref{lp}.
$\blacksquare$

By (\ref{lada}), there exists a pair $(w_{-4},\widetilde w_{-4})\in
H^1(\Omega)$ which solves  the boundary value problem
\begin{equation}
\mathcal M(x,D)(w_{-4},\widetilde w_{-4})=0\quad
\mbox{in}\,\,\Omega,\quad
(w_{-4}+\widetilde w_{-4})\vert_{\Gamma_0}=\frac{\widetilde
q_3}{2\partial_z\Phi}+\frac{\widetilde q_4}{2\partial_{\overline z}
\overline\Phi}.
\end{equation}

We set $(w_{-3,s},\widetilde w_{-3,s})=e^{s\eta}(w_{-3},\widetilde
w_{-3}), (\widetilde a_{\pm,s},\widetilde b_{\pm,s})
=e^{s\eta} (\widetilde a_{\pm},\widetilde b_{\pm}).$
We introduce the function $w_{1,s,\tau},
\widetilde w_{1,s,\tau}$ by formulas
\begin{equation}\label{-1}
w_{1,s,\tau}=w_{1,s}+\frac{w_{-3,s}+e_2q_{3,s}/2\partial_{ z}
\Phi}{\tau}+\frac{1}{\tau^2}( e^{2i\tau\psi(\widetilde x)}\widetilde
a_{+,s}+e^{-2i\tau\psi(\widetilde x)}\widetilde
a_{-,s}+w_{-4,s}-\frac{e_2\widetilde q_{3,s}}{2\partial_z\Phi})
\end{equation}
and
\begin{equation}\label{-2}
\widetilde w_{1,s,\tau}=\widetilde w_{1,s}+\frac{\widetilde
w_{-3,s}+e_2q_{4,s}/2\partial_{\overline z}\overline
\Phi}{\tau}+\frac{1}{\tau^2}( e^{2i\tau\psi(\widetilde x)}\widetilde
b_{+,s}+e^{-2i\tau\psi(\widetilde x)}\widetilde b_{-,s}+\widetilde
w_{-4,s}-\frac{e_2\widetilde q_{4,s}}{2\partial_{\overline z}\overline\Phi}).
\end{equation}
By (\ref{lada1}) and (\ref{lada}), the functions $w_{1,s,\tau}, \widetilde
w_{1,s,\tau}$ belong to $H^1(\Omega).$ Using (\ref{001}),  for any
$p\in (1,+\infty)$ we have
\begin{eqnarray}\label{90}
{L}_2(x,s,{D})^{*}\left(-e^{-\tau\Phi}\widetilde{\mathcal R}
_{-\tau,-A_{2,s}^*}(e_1(q_{3,s}+\frac{\widetilde
q_{3,s}}{\tau}))+\frac{e^{-\tau\Phi}e_2(q_{3,s}+\frac{\widetilde
q_{3,s}}{\tau})}{2\tau\partial_z\Phi}\right.\nonumber\\
\left.-e^{-\tau\overline\Phi}{\mathcal R}_{-\tau,- B_{2,s}^*}
(e_1(q_{4,s}+\frac{\widetilde
q_{4,s}}{\tau}))+\frac{e^{-\tau\overline\Phi}e_2(q_{4,s}+\frac{\widetilde
q_{4,s}}{\tau})}{2\tau\partial_{\overline z}\overline\Phi}\right )
\nonumber\\
= -{L}_2(x,s,{D})^{*}\left (e^{-\tau\Phi}\widetilde{\mathcal R}_{-\tau,-
A_{2,s}^*}(e_1(q_{3,s}+\frac{\widetilde
q_{3,s}}{\tau}))-\frac{e^{-\tau\Phi}e_1(q_{3,s}+\frac{\widetilde
q_{3,s}}{\tau})}{2\tau\partial_z\Phi}\right)                     \nonumber\\
- {L}_2(x,s,{D})^{*}\left (e^{-\tau\overline\Phi}{\mathcal
R}_{-\tau,-B_{2,s}^*}(e_1(q_{4,s}+\frac{\widetilde
q_{4,s}}{\tau}))-\frac{e^{-\tau\overline\Phi}e_2(q_{4,s}+\frac{\widetilde
q_{4,s}}{\tau})}{2\tau\partial_{\overline z}\overline\Phi}\right )
                                        \nonumber\\
= -e^{\tau\Phi}(Q_2^* -2 \partial_z {{B}^*_{2}} -{{B}_{2}}^*
{{A}_{2}}^*)\widetilde{\mathcal R}_{-\tau,-
A_{2,s}^*}(e_1(q_{3,s}+\frac{\widetilde
q_{3,s}}{\tau}))+e^{\tau\Phi}L_2(x,s,D)^*(\frac{e_2(q_{3,s}+\frac{\widetilde
q_{3,s}}{\tau})}{2\tau\partial_z\Phi}) \nonumber \\
-e^{-\tau\overline\Phi}(Q_2^* -2 \partial_{\overline z} {{A}^*_2} -{{A}_2}^*
{{B}_2}^*){\mathcal R}_{-\tau,-B_{2,s}^*}(e_1(q_{4,s}+\frac{\widetilde
q_{4,s}}{\tau}))
+e^{\tau\overline\Phi}L_2(x,s,D)^*(\frac{e_2(q_{4,s}+\frac{\widetilde
q_{4,s}}{\tau})}{2\tau\partial_{\overline z}\overline\Phi})\nonumber\\
-( Q_2^* -2\partial_{\overline z} {{A}_2}^* -{{A}_2}^* {{B}_2}^*
)(w_{1,s}+\frac{w_{-3,s}}{\tau})e^{-\tau\Phi}\nonumber\\ -( Q_2^* -2
\partial_z {{B}^*_2} -{{B}_2}^* {{A}_2}^*)(\widetilde
w_{1,s}+\frac{w_{-3,s}}{\tau})
e^{-\tau\overline\Phi}+o_{L^p(\Omega)}(\frac{1}{\tau}).
\end{eqnarray}

Setting $v^*=w_{1,s,\tau}e^{-\tau \Phi}+\widetilde w_{1,s,\tau}e^{-\tau
\overline \Phi}-e^{-\tau\Phi}\widetilde{\mathcal R}_{-\tau,-
A^*_{2,s}}(e_1(q_{3,s}+\frac{\widetilde
q_{3,s}}{\tau}))-e^{-\tau\overline\Phi}{\mathcal R}_{-\tau,
-B_{2,s}^*}(e_1(q_{4,s}+\frac{\widetilde q_{4,s}}{\tau}))$  for any
$p\in (1,\infty)$, we obtain that

\begin{equation}\label{nino}
L_2(x,s,D)v^*=e^{-\tau\varphi}o_{L^p(\Omega)}(\frac
1\tau)\quad\mbox{in}\,\,\Omega, \quad
v^*\vert_{\Gamma_0}=e^{-\tau\varphi}o_{H^1(\Gamma_0)}(\frac 1\tau).
\end{equation}

Using (\ref{nino}) and Proposition \ref{Proposition 0} and
\ref {Proposition 00}, we construct the last term $v_{-1}$ in
complex geometric optics solution which solves
the boundary value problem
\begin{equation}\label{nino1}
L_2(x,s,D)v_{-1}=L_2(x,s,D)v^* \quad\mbox{in}\,\,\Omega, \quad
v_{-1}\vert_{\Gamma_0}=v^* \end{equation}
and we obtain
\begin{equation}\label{Amimino11}
\root\of{\vert\tau\vert} \Vert v_{-1} \Vert_{L^2(\Omega)} +
\frac{1}{\root\of{\vert\tau\vert} } \Vert (\nabla v_{-1})
\Vert_{L^2(\Omega)}+\Vert v_{-1} \Vert_{H^{1,\tau}(\mathcal
O_{\epsilon})}=o(\frac{1}{\tau}).
\end{equation}
Finally we have a complex geometric optics solution for
Schr\"odinger operator $L_2(x,s,D)$ in a form:
\begin{equation}\label{-3}
v=w_{1,s,\tau}e^{-\tau \Phi}+\widetilde w_{1,s,\tau}e^{-\tau \overline
\Phi}-e^{-\tau\Phi}\widetilde{\mathcal R}_{-\tau,-
A^*_{2,s}}(e_1(q_{3,s}+\frac{\widetilde
q_{3,s}}{\tau}))
\end{equation}
$$
-e^{-\tau\overline\Phi}{\mathcal R}_{-\tau,
-B_{2,s}^*}(e_1(q_{4,s}+\frac{\widetilde q_{4,s}}{\tau}))
+v_{-1}e^{-\tau \varphi}.
$$

By (\ref{-3}), (\ref{nino}) and (\ref{nino1}), we have
\begin{equation}\label{-4}
L_2(x,s,D)v=0\quad\mbox{in}\,\,\Omega, \quad v\vert_{\Gamma_0}=0.
\end{equation}

{\bf Step 2:Asymptotic formula.}

Let $u_2=u_2(s,x)$ be a solution to the following boundary value
problem:
\begin{equation}\label{(2.1I)}
{ L}_{2}(x,s,D)u_2=0\quad \mbox{in}\,\,\Omega,\quad
u_2\vert_{\partial\Omega}=u_1\vert_{\partial\Omega}, \quad
\frac{\partial u_2}{\partial \nu}\vert_{\widetilde \Gamma}
=\frac{\partial u_1}{\partial\nu}\vert_{\widetilde \Gamma}.
\end{equation}

Setting $u=u_1-u_2$, we have
\begin{equation}
{L}_2(x,s,{D})u +2({A}_1-{A}_{2,s})\partial_z u_1
+2({B}_1-{B}_{2,s})\partial_{\overline z} u_1
+(Q_1-Q_{2,s})u_1=0 \quad \mbox{in}~ \Omega    \label{mn}
\end{equation}
and
\begin{equation}\label{mn1}
u \vert_{\partial\Omega} =0, \quad \frac{\partial u}{\partial \nu}
\vert_{\widetilde \Gamma} =0.
\end{equation}
Let $v$ be a function given by  (\ref{-3}).  Taking the scalar
product of (\ref{mn}) with $ v$  in $L^2(\Omega)$ and using
(\ref{-4}) and (\ref{mn1}), we obtain
\begin{equation}\label{ippolit}
0=\frak G(u_1,v)= \int_{\Omega}(2({A}_1-{A}_{2,s})\partial_z
u_1 +2({B}_1-{B}_{2,s})\partial_{\overline z} u_1 +(Q_1-Q_{2,s})u_1,
\overline v) dx.
\end{equation}
Our goal is to obtain the asymptotic formula for the right-hand side
of (\ref{ippolit}). We have
\begin{proposition}\label{Nova}
There exists a constant $\mathcal C_0$, independent of $\tau$,
such that the following asymptotic formula is valid as $\vert
\tau\vert\rightarrow +\infty$:
\begin{eqnarray}\label{nonsence1}
I_0=((Q_1-Q_{2,s})u_1,v)_{L^2(\Omega)}\\
= \int_\Omega
(((Q_1-Q_{2,s})w_0,\overline
{\widetilde w_{1,s}})+((Q_1-Q_{2,s})\widetilde w_0,\overline{ w_{1,s}}) )dx
                                             \nonumber\\
+ \frac{\mathcal C_0}{\tau}
+2\pi \frac{((Q_1-Q_{2,s})w_0,\overline{ w_{1,s}})(\widetilde x) e^{2\tau
i\psi(\widetilde x)} +((Q_1-Q_{2,s})\widetilde w_0,\overline {\widetilde
w_{1,s}})(\widetilde x)e^{-2i\tau\psi(\widetilde x)}} {\tau\vert
\mbox{det}\, \psi ''(\widetilde x)\vert^\frac 12}
                                               \nonumber\\
+\frac {1}{2\tau i}\int_{\partial\Omega} ((Q_1-Q_{2,s})w_0,
\overline{w_{1,s}})e^{2\tau i\psi}
\frac{(\nu,\nabla\psi)}{\vert\nabla\psi\vert^2}d\sigma - \frac
{1}{2\tau i}\int_{\partial\Omega} ((Q_1-Q_{2,s})\widetilde w_0,\overline
{\widetilde w_{1,s}}) e^{-2\tau i\psi}
\frac{(\nu,\nabla\psi)}{\vert\nabla\psi\vert^2}d\sigma
+ o(\frac 1\tau).
                                       \nonumber
\end{eqnarray}
\end{proposition}
{\bf Proof.}
By (\ref{zad1}), (\ref{zad2}), (\ref{kl}), (\ref{gandon1}), (\ref{zad}) and
Propositions \ref{vasya} and \ref{gandon}, we have
\begin{equation}\label{P1}
u_1(x)=(w_0+\frac{w_{-1}}{\tau})e^{\tau\Phi}+({\widetilde
w_0}+\frac{\widetilde
w_{-1}}{\tau})e^{\tau\overline\Phi}-\frac{q_2e^{\tau\overline\Phi}}
{2\tau\overline
{\partial_z\Phi}}-\frac{q_1e^{\tau\Phi}}
{2\tau{\partial_z\Phi}}+e^{\tau\varphi}o_{L^2(\Omega)}
(\frac 1\tau)\,\,\mbox{as}\,\,\tau\rightarrow +\infty.
\end{equation}
Using (\ref{-1}), (\ref{-2}), (\ref{lada1}), (\ref{lada}), (\ref{-3})
and Propositions \ref{vasya} and \ref{gandon}, we obtain
\begin{equation}\label{P2}
v(x)=(w_{1,s}+\frac{w_{-2,s}}{\tau})e^{-\tau\Phi}+(\widetilde
w_{1,s}+\frac{\widetilde w_{-2,s}}{\tau})
e^{-\tau\overline\Phi}+\frac{q_{4,s}e^{-\tau\overline\Phi}}
{2\tau\overline{\partial_z\Phi}}+\frac{q_{3,s}e^{-\tau\Phi}}
{2\tau{\partial_z\Phi}}+e^{-\tau\varphi}o_{L^2(\Omega)}(\frac
1\tau)\,\,\mbox{as}\,\,\tau\rightarrow +\infty.
\end{equation}

By (\ref{P1}) and (\ref{P2}), we obtain the following asymptotic formula:
\begin{eqnarray}\label{nonsence}
((Q_1-Q_{2,s})u_1,v)_{L^2(\Omega)} =
((Q_1-Q_{2,s})((w_0+\frac{w_{-1}}{\tau})e^{\tau\Phi}+({\widetilde
w_0}+\frac{\widetilde
w_{-1}}{\tau})e^{\tau\overline\Phi}-\frac{q_2e^{\tau\overline\Phi}}
{2\tau\overline{\partial_z\Phi}}-\frac{q_1e^{\tau\Phi}}
{2\tau{\partial_z\Phi}}+e^{\tau\varphi}o_{L^2(\Omega)}(\frac 1\tau)),
\nonumber\\
(w_{1,s}+\frac{w_{-2,s}}{\tau})e^{-\tau\Phi}+({\widetilde
w_{1,s}}+\frac{\widetilde
w_{-2,s}}{\tau})e^{-\tau\overline\Phi}+\frac{q_{4,s}e^{-\tau\overline\Phi}}
{2\tau\overline{\partial_z\Phi}}+\frac{q_{3,s}e^{-\tau\Phi}}
{2\tau{\partial_z\Phi}}+e^{-\tau\varphi}o_{L^2(\Omega)}
(\frac 1\tau))
_{L^2(\Omega)}                              \nonumber\\
= \int_\Omega ( ((Q_1-Q_{2,s}) \widetilde w_0,\overline{w_{1,s}})
+ \frac 1\tau((Q_1-Q_{2,s})\widetilde w_0, \overline{w_{-2,s}})
+ \frac
{1}{\tau}((Q_1-Q_{2,s})\widetilde w_{-1}, \overline{w_{1,s}})
                                 \nonumber\\
+((Q_1-Q_{2,s}) w_0,\overline{\widetilde w_{1,s}})
+\frac{1}{\tau}((Q_1-Q_{2,s})w_{-1},\overline{\widetilde
w_{1,s}})+((Q_1-Q_{2,s})w_0, \overline{w_{-2,s}}))dx\nonumber\\
+\frac 1\tau\int_\Omega ((Q_1-Q_{2,s})w_0,\frac{\overline{q_{4,s}}
}{2{\partial_z\Phi}})-((Q_1-Q_{2,s})\frac{{q_2}}{2\partial_z\Phi},
\overline{w_{1,s}}) \nonumber \\
- ((Q_1-Q_{2,s})\frac{{q_1}} {2\partial_{\overline z}
\overline\Phi}, \overline{\widetilde w_{1,s}})
+ ((Q_1-Q_{2,s})\widetilde w_0,\frac{\overline{q_{3,s}}
}{2\overline{\partial_z\Phi}}))dx\nonumber\\
+ \int_\Omega (((Q_1-Q_{2,s})w_0,\overline{ w_{1,s}})
e^{2i\tau\psi}
+ ((Q_1-Q_{2,s})\widetilde w_0, \overline{\widetilde w_{1,s}})
e^{-2i\tau\psi})dx +o(\frac{1}{\tau}).\nonumber
\end{eqnarray}
Applying the stationary phase argument (see e.g., \cite{BH}) to the
last integral on the right-hand side of this formula, we complete the
proof of Proposition \ref{Nova}. $\blacksquare$

We set
$$
\mathcal{U}=w_{0,\tau}e^{ \tau \Phi} +\widetilde w_{0,\tau}e^{\tau
\overline{\Phi}}, \,\, \mathcal{V}=w_{1,s,\tau}e^{- \tau \Phi}
+\widetilde w_{1,s,\tau} e^{-\tau \overline{\Phi}}.
$$
By the stationary phase argument and formulae (\ref{-5}),
(\ref{ll1}), (\ref{zad1}), (\ref{zad2}), (\ref{-1}) and (\ref{-2}),
short calculations yield that there exist constants
$\kappa_k,\widetilde \kappa_k$, independent of $\tau$, such that
\begin{eqnarray}\label{MAHA}
I_1 \equiv 2(({A}_1-{A}_{2,s})\partial_z \mathcal{U},\mathcal{V})
_{L^2(\Omega)} \nonumber\\
=  (2({A}_1-{A}_{2,s})( \partial_z (w_{0,\tau}e^{\tau \Phi})
+\partial_z {\widetilde w_{0,\tau}}e^{\tau \overline{\Phi}}),
 w_{1,s,\tau} e^{-\tau \Phi} +\widetilde w_{1,s,\tau}
e^{-\tau \overline{\Phi}})_{L^2(\Omega)}                           \nonumber\\
=\sum_{k=1}^3\tau^{2-k}\kappa_k+ e^{2i\tau\psi(\widetilde
x)}(((A_1-A_{2,s})\partial_z\Phi a_+,\widetilde
w_{1,s})_{L^2(\Omega)}+((A_1-A_{2,s})\partial_z\Phi w_0,\widetilde
b_{-,s})_{L^2(\Omega)})\nonumber\\+ e^{-2i\tau\psi(\widetilde
x)}(((A_1-A_{2,s})\partial_z\Phi a_-, \widetilde
w_{1,s})_{L^2(\Omega)}+((A_1-A_{2,s})\partial_z\Phi w_0, \widetilde
b_{+,s})_{L^2(\Omega)})\nonumber\\
+2\int_{\Omega} (({A}_1-{A}_{2,s}) \partial_z {\widetilde
w_0},\overline {{\widetilde w}_{1,s}})e^{-2i \tau \psi}dx  -
\int_\Omega (2\partial_z({A}_1-{A}_{2,s}) w_0, \overline{
w_{1,s}})e^{2\tau i\psi}dx
\nonumber\\
-  \int_\Omega 2(({A}_1-{A}_{2,s})w_0 ,\overline{\partial_{\overline
z}w_{1,s}})e^{2\tau i\psi}dx
+ \int_{\partial\Omega}(\nu_1-i\nu_2)((A_1-A_{2,s})w_0,\overline{w_{1,s}})
e^{2i\tau\psi}d\sigma \nonumber\\
+ o\left( \frac{1}{\tau}\right) \nonumber\\
=\sum_{k=1}^3\tau^{2-k}\kappa_k+ e^{2i\tau\psi(\widetilde
x)}(((A_1-A_{2,s})\partial_z\Phi a_+,\widetilde
w_{1,s})_{L^2(\Omega)}+((A_1-A_{2,s})\partial_z\Phi w_0,\widetilde
b_{-,s})_{L^2(\Omega)})\nonumber\\+ e^{-2i\tau\psi(\widetilde
x)}(((A_1-A_{2,s})\partial_z\Phi a_-,  \widetilde
w_{1,s})_{L^2(\Omega)}+((A_1-A_{2,s})\partial_z\Phi w_0, \widetilde
b_{+,s})_{L^2(\Omega)})\nonumber\\ -\int_{\Omega} (({A}_1-{A}_{2,s})
B_1\widetilde w_0,\overline {{\widetilde w}_{1,s}})e^{-2i \tau
\psi}dx  - \int_\Omega (2{\partial_ z}({A}_1-{A}_{2,s}) w_0,
\overline{ w_{1,s}})e^{2\tau i\psi}dx
\nonumber\\
-  \int_\Omega (({A}_1-{A}_{2,s})w_0
,\overline{B_{2,s}^*w_{1,s}})e^{2\tau i\psi}dx
+\int_{\partial\Omega}(\nu_1-i\nu_2)((A_1-A_{2,s})w_0,\overline{w_{1,s}})
e^{2i\tau\psi}d\sigma \nonumber \\
+ o\left( \frac{1}{\tau}\right)
\end{eqnarray}
\nopagebreak and
\begin{eqnarray}\label{MAHA1}
I_2 \equiv ( ({B}_1-{B}_{2,s}){\partial_{\overline z}
\mathcal{U}},\mathcal{V})
_{L^2(\Omega)}                                   \nonumber\\
= (2({B}_1-{B}_{2,s})(e^{\tau \Phi} \partial_{\overline z} w_{0,\tau} +
{\partial}_{ \overline{z}} (\widetilde w_{0,\tau}e^{\tau
\overline{\Phi}})),
  w_{1,s,\tau} e^{ -\tau \Phi}
+ \widetilde w_{1,s,\tau} e^{-\tau \overline{\Phi}})_{L^2(\Omega)}
                                   \nonumber \\
=\sum_{k=1}^3\tau^{2-k}\widetilde\kappa_k + e^{2i\tau\psi(\widetilde
x)}(((B_1-B_{2,s})\partial_{\overline z}\overline\Phi b_+,
w_{1,s})_{L^2(\Omega)}+((B_1-B_{2,s})\partial_{\overline z}\overline\Phi
\widetilde w_0,\widetilde a_{-,s})_{L^2(\Omega)})
\nonumber\\
+ e^{-2i\tau\psi(\widetilde
x)}(((B_1-B_{2,s})\partial_{\overline z}\overline\Phi b_-,
w_{1,s})_{L^2(\Omega)}+((B_1-B_{2,s})\partial_{\overline z}\overline\Phi
\widetilde w_0, \widetilde a_{+,s})_{L^2(\Omega)})
\nonumber\\
+ \int_{\Omega}2(({B}_1-{B}_{2,s}) {\partial_{\overline z} w_0},
\overline{w_{1,s}}) e^{2
\tau i \psi}dx -\int_\Omega(2\partial_{\overline{z}}
({B}_1-{B}_{2,s}) \widetilde w_0, \overline{\widetilde w_{1,s}}) e^{-2\tau i
\psi}dx                              \nonumber\\
-\int_\Omega(2({B}_1-{B}_{2,s})\widetilde w_0,
 \overline{\partial_z\widetilde w_{1,s}}) e^{-2\tau \psi}dx
+\int_{\partial\Omega}(\nu_1+i\nu_2)((B_1-B_{2,s}) \widetilde w_0,
\overline{\widetilde w_{1,s}})e^{ -2i\tau\psi}d\sigma
+ o\left(\frac{1}{\tau}\right)
                                                   \nonumber\\
=\sum_{k=1}^3\tau^{2-k}\widetilde\kappa_k + e^{2i\tau\psi(\widetilde
x)}(((B_1-B_{2,s})\partial_{\overline z}\overline\Phi b_+,
w_{1,s})_{L^2(\Omega)}+((B_1-B_{2,s})\partial_{\overline z}\overline\Phi
\widetilde w_0, a_{-,s})_{L^2(\Omega)})\nonumber\\+ e^{-2i\tau\psi(\widetilde
x)}(((B_1-B_{2,s})\partial_{\overline z}\overline\Phi b_-,
w_{1,s})_{L^2(\Omega)}+((B_1-B_{2,s})\partial_{\overline z}\overline\Phi
\widetilde w_0, a_{+,s})_{L^2(\Omega)})\nonumber\\
- \int_{\Omega}
(({B}_1-{B}_{2,s}) A_1w_0, \overline{w_{1,s}}) e^{2 \tau i \psi}dx
-\int_\Omega(2{\partial}_{ \overline{z}} ({B}_1-{B}_{2,s}) \widetilde
w_0, \overline{\widetilde w_{1,s}}) e^{-2\tau i \psi}dx
\nonumber\\
-\int_\Omega(({B}_1-{B}_{2,s})\widetilde w_0,\overline{
A_{2,s}^* \widetilde w_{1,s}}) e^{-2\tau i \psi}dx
+\int_{\partial\Omega}(\nu_1+i\nu_2)((B_1-B_{2,s}) \widetilde w_0,
\overline{\widetilde w_{1,s}})e^{ -2i\tau\psi}d\sigma + o\left(
\frac{1}{\tau}\right).
\end{eqnarray}

Using (\ref{0011}) and integrating by parts, we obtain
\begin{eqnarray}\label{-11}
I_3=-\int_\Omega (2(A_1-A_{2,s})\partial_z(e^{\tau\Phi}\widetilde
{\mathcal R}_{\tau,B_1} \{e_1(q_1+\widetilde q_1/\tau)\}+e^{\tau
\overline\Phi} {\mathcal R}_{\tau,A_1}
\{e_1(q_2+\widetilde q_2/\tau )\})
                                  \nonumber\\
+ 2(B_1-B_{2,s})\partial_{\overline z}(e^{\tau\Phi}\widetilde
{\mathcal R}_{\tau,B_1} \{e_1(q_1+\widetilde q_1/\tau)\}+e^{\tau
\overline\Phi} {\mathcal R}_{\tau,A_1} \{e_1(q_2+\widetilde q_2/\tau )\}),
\overline{\mathcal V})dx                 \nonumber \\
= -\int_\Omega (2(A_1-A_{2,s})(e^{\tau\Phi}(-B_1\widetilde {\mathcal
R}_{\tau,B_1} \{e_1(q_1+\widetilde q_1/\tau)\}+e_1(q_1+\widetilde q_1/\tau
))+e^{\tau \overline\Phi}
\partial_z{\mathcal
R}_{\tau,A_1} \{e_1(q_2+\widetilde q_2/\tau )\})\nonumber\\
+ 2(B_1-B_{2,s})(e^{\tau\Phi}\partial_{\overline z}\widetilde {\mathcal
R}_{\tau,B_1} \{e_1(q_1+\widetilde q_1/\tau )\}+e^{\tau \overline\Phi}(-A_1
{\mathcal R}_{\tau,A_1}
\{e_1(q_2+\widetilde q_2/\tau )\}+e_1(q_2+\widetilde q_2/\tau ))),
\overline {\mathcal V})dx\nonumber\\
= -\int_\Omega (2(A_1-A_{2,s})e^{\tau\Phi}(-B_1\widetilde {\mathcal
R}_{\tau,B_1} \{e_1(q_1+\widetilde q_1/\tau)\}+e_1(q_1+\widetilde
q_1/\tau))\nonumber\\
+2(B_1-B_{2,s})e^{\tau \overline\Phi}(-A_1 {\mathcal
R}_{\tau,A_1} \{e_1(q_2+\widetilde q_2/\tau)\}+e_1(q_2+\widetilde q_2/\tau)),
\overline {\mathcal V})dx           \nonumber\\
+ \int_\Omega (2\partial_z(A_1-A_{2,s}) e^{\tau\overline\Phi}\mathcal
R_{\tau,A_1}\{e_1(q_2+\widetilde q_2/\tau)\}+2\partial_{\overline
z}(B_1-B_{2,s})e^{\tau\Phi}\widetilde{\mathcal R}_{\tau,
B_1}\{e_1(q_1+\widetilde q_1/\tau)\},
\overline {\mathcal V})dx\nonumber\\
+\int_\Omega (2(A_1-A_{2,s}) e^{\tau\overline\Phi}\mathcal
R_{\tau,A_1}\{e_1(q_2+\widetilde q_2/\tau)\},
\partial_z \overline {\mathcal V})+(2(B_1-B_{2,s})e^{\tau\Phi}
\widetilde{\mathcal R}_{\tau,
B_1}\{e_1(q_1+\widetilde q_1/\tau)\},
\partial_{\overline z}\overline {\mathcal V})dx\nonumber\\
-\int_{\partial\Omega} \{(\nu_1-i\nu_2)((A_1-A_{2,s})
e^{\tau\overline\Phi}\mathcal R_{\tau,A_1}\{e_1(q_2+\widetilde q_2/\tau )\},
\overline {\mathcal V})
                                             \nonumber\\
+(\nu_1+i\nu_2)((B_1-B_{2,s})e^{\tau\Phi}\widetilde{\mathcal
R}_{\tau, B_1}\{e_1(q_1+\widetilde q_1/\tau)\},
\overline {\mathcal V})\}d\sigma \nonumber\\
= -\int_\Omega (2(A_1-A_{2,s})e^{\tau\Phi}(-B_1\widetilde {\mathcal
R}_{\tau,B_1} \{e_1(q_1+\widetilde q_1/\tau )\}+e_1(q_1+\widetilde
q_1/\tau ))\nonumber\\+2(B_1-B_{2,s})e^{\tau \overline\Phi}(-A_1
{\mathcal
R}_{\tau,A_1} \{ e_1(q_2+\widetilde q_2/\tau )\}
+e_1(q_2+\widetilde q_2/\tau )), \overline {\mathcal V})dx
                             \nonumber\\
+\int_\Omega (2\partial_z(A_1-A_{2,s}) e^{\tau\overline\Phi}\mathcal
R_{\tau,A_1}\{e_1(q_2+\widetilde q_2/\tau)\}+2\partial_{\overline
z}(B_1-B_{2,s})e^{\tau\Phi}\widetilde{\mathcal R}_{\tau,
B_1}\{e_1(q_1+\widetilde q_1/\tau)\},
\overline {\mathcal V})dx\nonumber\\
+\int_\Omega (2(A_1-A_{2,s}) \mathcal R_{\tau,A_1}\{e_1(q_2+\widetilde
q_2/\tau)\},
\partial_z \overline {w_{1,s}})+(2(B_1-B_{2,s})\widetilde{\mathcal R}_{\tau,
B_1}\{e_1(q_1+\widetilde q_1/\tau )\},
\partial_{\overline z}\overline {\widetilde w_{1,s}})dx\nonumber\\
+2\int_\Omega e^{\tau (\overline \Phi-\Phi)}(e_1(q_2+\widetilde q_2/\tau
),\overline{
\mbox{\bf P}^*_{A_1}((A_1-A_{2,s})^*
(\partial_{\overline z} {\widetilde w}_{1,s}
-\tau\partial_{\overline z}\overline\Phi{\widetilde w}_{1,s}))})dx
                           \nonumber\\
+2\int_\Omega e^{\tau ( \Phi-\overline \Phi)}(e_1(q_1+\widetilde q_1/\tau ),
 \overline{\mbox{\bf T}^*_{B_1}((B_1-B_{2,s})^*(\partial_{ z} w_{1,s}
-\tau {\partial_z\Phi} w_{1,s}))})dx
                                    \nonumber\\
-\int_{\partial\Omega}\{ (\nu_1-i\nu_2)((A_1-A_{2,s})
e^{\tau\overline\Phi}\mathcal R_{\tau,A_1}\{e_1(q_2+\widetilde
q_2/\tau )\}, \overline{\mathcal V})
                                \nonumber\\
+(\nu_1+i\nu_2)((B_1-B_{2,s})e^{\tau\Phi}\widetilde{\mathcal
R}_{\tau, B_1}\{e_1(q_1+\widetilde q_1/\tau)\}, \overline {\mathcal
V})\}d\sigma.
\end{eqnarray}

By (\ref{50l}) and (\ref{50ll}), the boundary integrals in (\ref{-11})
are $O(\frac{1}{\tau^2}).$
By (\ref{gandon1}) and Proposition \ref{gandon}, we have
\begin{eqnarray}\label{dd}
2\int_\Omega e^{\tau (\overline \Phi-\Phi)}(e_1\widetilde q_2/\tau
,\overline{
\mbox{\bf P}^*_{A_1}((A_1-A_{2,s})^*(\partial_{\overline z}
{\widetilde w}_{1,s}-\tau\partial_{\overline z}
\overline\Phi{\widetilde w}_{1,s}))})dx\nonumber\\
+ 2\int_\Omega e^{\tau ( \Phi-\overline \Phi)}(e_1\widetilde q_1/\tau ,
 \overline{\mbox{\bf T}^*_{B_1}((B_1-B_{2,s})^*(\partial_{ z} w_{1,s}
-\tau {\partial_z\Phi} w_{1,s}))})dx=o(\frac{1}{\tau})\quad
\mbox{as}\,\,\tau\rightarrow +\infty.
\end{eqnarray}
Applying the stationary phase argument, (\ref{dd}),  Propositions
\ref{vasya} and \ref{gandonnal}, we obtain from (\ref{-11}) that
there exists a constant $\mathcal C_1$ independent of $\tau$ such
that
\begin{eqnarray}
I_3=\frac{\mathcal C_1}{\tau}+2\int_\Omega e^{\tau (\overline
\Phi-\Phi)}(e_1q_2,\overline{
\mbox{\bf P}^*_{A_1}((A_1-A_{2,s})^*(-\tau\partial_{\overline z}
\overline\Phi{\widetilde w}_{1,s}))})dx\nonumber\\
+2\int_\Omega e^{\tau ( \Phi-\overline \Phi)}(e_1q_1,
 \overline{\mbox{\bf T}^*_{B_1}((B_1-B_{2,s})^*(-\tau {\partial_z\Phi}
w_{1,s}))})dx+
o(\frac{1}{\tau})\quad \mbox{as}\,\,\tau\rightarrow +\infty.
\end{eqnarray}

Using (\ref{0011}) and integrating by parts, we obtain
\begin{eqnarray}\label{-10}
I_4=\int_\Omega (2(A_1-A_{2,s})\partial_z\mathcal
U+2(B_1-B_{2,s})\partial_{\overline z} \mathcal
U,\nonumber\\\overline{-e^{-\tau\Phi}\widetilde{\mathcal R}_{-\tau,-
A_{2,s}^*}\{e_1(q_{3,s}+\widetilde q_{3,s}/\tau
)\}-e^{-\tau\overline\Phi}{\mathcal R}_{-\tau,-
B_{2,s}^*}\{e_1(q_{4,s}+\widetilde q_{4,s}/\tau )\}})dx
                                            \nonumber \\
= -\int_\Omega (2(A_1-A_{2,s})\partial_{ z}\widetilde w_0
e^{\tau\overline\Phi}+2(B_1-B_{2,s})\partial_{\overline z}
w_0 e^{\tau\Phi},\nonumber\\
\overline{e^{-\tau\Phi}\widetilde{\mathcal R}_{-\tau,-
A_2^*}\{e_1(q_{3,s}+\widetilde q_{3,s}/\tau
)\}+e^{-\tau\overline\Phi}{\mathcal R}_{-\tau,-
B_{2,s}^*}\{e_1(q_{4,s}+\widetilde q_{4,s}/\tau )\}})dx\nonumber\\
-\int_\Omega (2(A_1-A_{2,s})(\partial_z w_0+\tau\partial_z\Phi
w_0)e^{\tau\Phi}, \overline{e^{-\tau\Phi}\widetilde{\mathcal R}
_{-\tau,-A_{2,s}^*}\{e_1(q_{3,s}+\widetilde q_{3,s}/\tau
)\})}dx\nonumber\\-\int_\Omega(2(B_1-B_{2,s})(\partial_{\overline z}
\widetilde w_0+\tau\partial_{\overline z}\overline \Phi \widetilde
w_0)e^{\tau\overline\Phi}, \overline {e^{-\tau\overline\Phi}{\mathcal
R}_{-\tau,-
B_{2,s}^*}\{e_1(q_{4,s}+\widetilde q_{4,s}/\tau )\}})dx\nonumber\\
+\int_\Omega ((2\partial_z(A_1-A_{2,s})w_0
e^{\tau\Phi},\overline{e^{-\tau\overline\Phi}{\mathcal R}_{-\tau,-
B_2^*}\{e_1(q_{4,s}+\widetilde q_{4,s}/\tau
)\}})\nonumber\\+(2\partial_{\overline z}(B_1-B_{2,s})\widetilde w_0
e^{\tau\overline\Phi},\overline{e^{-\tau\Phi}\widetilde{\mathcal
R}_{-\tau,- A_{2,s}^*}\{e_1(q_{3,s}+\widetilde q_{3,s}/\tau )\}}))dx\nonumber\\
-\int_{\partial\Omega} \{(\nu_1-i\nu_2)((A_1-A_{2,s})w_0
e^{\tau\Phi},\overline{e^{-\tau\overline\Phi}{\mathcal R}_{-\tau,-
B_{2,s}^*}\{ e_1(q_{4,s}+\widetilde
q_{4,s}/\tau )\}})\nonumber\\
+(\nu_1+i\nu_2)((B_1-B_{2,s})\widetilde w_0
e^{\tau\overline\Phi},\overline{e^{-\tau\Phi}\widetilde{\mathcal
R}_{-\tau,-A_{2,s}^*}\{e_1(q_{3,s}+\widetilde q_{3,s}/\tau )\}})\}d\sigma
                                     \nonumber\\
+\int_\Omega (2(A_1-A_{2,s})w_0
e^{\tau\Phi},\overline{\partial_{\overline z}(e^{-\tau\overline\Phi}{\mathcal
R}_{-\tau,-B_{2,s}^*}\{e_1(q_{4,s}+\widetilde q_{4,s}/\tau )\})})dx
\nonumber\\
+\int_\Omega (2(B_1-B_{2,s})\widetilde w_0
e^{\tau\overline\Phi},\overline{\partial_z(e^{-\tau\Phi}\widetilde{\mathcal
R}_{-\tau,-A_{2,s}^*}\{e_1(q_{3,s}+\widetilde q_{3,s}/\tau )\}}))dx
                                               \nonumber\\
= -\int_\Omega (2(A_1-A_{2,s})\partial_{ z}\widetilde w_0
e^{\tau\overline\Phi}+2(B_1-B_{2,s})\partial_{\overline z} w_0 e^{\tau
\Phi}, \nonumber\\
\overline{e^{-\tau\Phi}\widetilde{\mathcal R}_{-\tau,-
A_{2,s}^*}\{e_1(q_{3,s}+\widetilde q_{3,s}/\tau
)\}+e^{-\tau\overline\Phi}{\mathcal
R}_{-\tau,- B_{2,s}^*}\{e_1(q_{4,s}+\widetilde q_{4,s}/\tau )\}})dx\nonumber\\
-\int_\Omega (2\mbox{\bf T}^*_{- A_{2,s}^*}((A_1-A_{2,s})(\partial_z
w_0+\tau\partial_z\Phi w_0)),
\overline{e^{\tau(\overline\Phi-\Phi)}e_1(q_{3,s}+\widetilde q_{3,s}/\tau
))}dx\nonumber\\-\int_\Omega(2\mbox{\bf P}^*_{-B_{2,s}^*}((B_1-B_{2,s})(\partial_{\overline
z} \widetilde w_0+\tau\partial_{\overline z}\overline \Phi \widetilde w_0)),
\overline {e^{\tau(\Phi-\overline\Phi)}e_1(q_{4,s}+\widetilde q_{4,s}/\tau )})dx\nonumber\\
+\int_\Omega (2\partial_z(A_1-A_{2,s})w_0
e^{\tau\Phi},\overline{e^{-\tau\overline\Phi}{\mathcal R}_{-\tau,-
B_{2,s}^*}\{e_1(q_{4,s}+\widetilde q_{4,s}/\tau )\}})\nonumber\\
+(2\partial_{\overline
z}(B_1-B_{2,s})\widetilde w_0
e^{\tau\overline\Phi},\overline{e^{-\tau\Phi}\widetilde{\mathcal
R}_{-\tau,- A_{2,s}^*}\{e_1(q_{3,s}+\widetilde q_{3,s}/\tau )\}}))dx\nonumber\\
-\int_{\partial\Omega} \{(\nu_1-i\nu_2)((A_1-A_{2,s})w_0
e^{\tau\Phi},\overline{e^{-\tau\overline\Phi}{\mathcal R}_{-\tau,-
B_{2,s}^*}\{e_1(q_{4,s}+\widetilde
q_{4,s}/\tau )\}})\nonumber\\
+(\nu_1+i\nu_2)((B_1-B_{2,s})\widetilde w_0
e^{\tau\overline\Phi},\overline{e^{-\tau\Phi}\widetilde{\mathcal
R}_{-\tau,- A_{2,s}^*}\{e_1(q_{3,s}+\widetilde q_{3,s}/\tau )\}})\}d\sigma
                                                 \nonumber\\
+\int_\Omega (2(A_1-A_{2,s})w_0 e^{\tau\Phi},\overline{
e^{-\tau\overline\Phi}B_{2,s}^*{\mathcal R}_{-\tau,-
B_{2,s}^*}\{e_1(q_{4,s}+\widetilde q_{4,s}/\tau )\}
+ e^{-\tau\overline\Phi}e_1(q_{4,s}+\widetilde q_{4,s}/\tau )})dx\nonumber\\
+\int_\Omega (2(B_1-B_{2,s})\widetilde w_0
e^{\tau\overline\Phi},\overline{e^{-\tau\Phi} A_{2,s}^*\widetilde{\mathcal
R}_{-\tau,- A_{2,s}^*}\{e_1(q_{3,s}+\widetilde q_{3,s}/\tau )\}
+e^{-\tau\Phi}e_1(q_{3,s}+\widetilde q_{3,s}/\tau )})dx.
\end{eqnarray}
By Proposition \ref{lp}, the boundary integral in (\ref{-10}) is
$O(\frac {1}{\tau^2}).$
By (\ref{lada}) and Proposition \ref{gandon}, we have
\begin{eqnarray}\label{-1111}
\frac 1\tau \int_\Omega (2\mbox{\bf T}^*_{- A_{2,s}^*}((A_1-A_{2,s})
(\partial_z
w_0+\tau\partial_z\Phi w_0)),
\overline{e^{\tau(\overline\Phi-\Phi)}e_1\widetilde
q_{3,s})}dx\nonumber\\-\int_\Omega(2\mbox{\bf P}^*_{-B_{2,s}^*}
((B_1-B_{2,s})(\partial_{\overline z} \widetilde w_0
+\tau\partial_{\overline z}\overline \Phi \widetilde w_0)), \overline
{e^{\tau(\Phi-\overline\Phi)}e_1\widetilde q_{4,s}})dx=o(\frac{1}{\tau})\quad
\mbox{as}\,\,\tau\rightarrow +\infty.
\end{eqnarray}

Applying the stationary phase argument,  Propositions \ref{vasya} and
\ref{gandonnal}, and (\ref{-1111}), we obtain  from (\ref{-10})
that there exists a constant $\mathcal C_2$, independent of $\tau$, such that
\begin{eqnarray}
I_4=\frac{\mathcal C_2}{\tau}-\int_\Omega (2\mbox{\bf T}^*_{-
A_{2,s}^*}((A_1-A_{2,s})\tau\partial_z\Phi w_0),
\overline{e^{\tau(\overline\Phi-\Phi)}e_1q_{3,s})}dx\nonumber\\
- \int_\Omega(2\mbox{\bf P}^*_{-B_{2,s}^*}
((B_1-B_{2,s})\tau\partial_{\overline z}\overline \Phi \widetilde
w_0), \overline
{e^{\tau(\Phi-\overline\Phi)}e_1q_{4,s}})dx+o(\frac{1}{\tau})\quad
\mbox{as}\,\,\tau\rightarrow +\infty.
\end{eqnarray}

{\bf Step 3:derivation of equations (\ref{op!})-(\ref{A2}).}

We set
\begin{eqnarray}\label{narkotic}
\mathcal{U}_1(x)=w_{0,\tau}e^{ \tau \Phi} +\widetilde w_{0,\tau}e^{\tau
\overline{\Phi}}-e^{\tau\Phi}\widetilde{\mathcal
R}_{\tau, B_1}\{e_1(q_1+\widetilde q_1/\tau )\}-e^{\tau\overline\Phi}{\mathcal R}_{\tau, A_1}\{ e_1(q_2+\widetilde q_2/\tau )\}, \nonumber\\
\mathcal{V}_1(x)=w_{1,s,\tau}e^{- \tau \Phi}
+ \widetilde w_{1,s,\tau}e^{-\tau \overline{\Phi}}
- e^{-\tau\Phi}\widetilde{\mathcal R}_{-\tau,-
A_{2,s}^*}\{e_1(q_{3,s}+\widetilde q_{3,s}/\tau )\}\nonumber\\
- e^{-\tau\overline \Phi}{\mathcal R}_{-\tau,- B_{2,s}^*}
\{ e_1(q_{4,s}+\widetilde
q_{4,s}/\tau )\}.\nonumber
\end{eqnarray}

By (\ref{mimino11}), (\ref{-3}) and Proposition \ref{vasya}, we have
\begin{equation}
\frak G(u_{-1}e^{\tau\varphi},v-(w_{1,s,\tau}e^{- \tau \Phi}
+ \widetilde w_{1,s,\tau}
e^{-\tau \overline{\Phi}}))=\frak G(u_1-(w_{0,\tau}e^{ \tau \Phi}
+\widetilde w_{0,\tau} e^{\tau
\overline{\Phi}}),v_{-1}e^{-\tau\varphi})
\end{equation}
$$
= o(\frac{1}{\root\of{\tau}})\quad\mbox{as}\,\,\tau\rightarrow
+\infty.
$$

Then
\begin{equation}
\frak G(u_1,v)=\int_{\widetilde\Gamma}(\nu_1-i\nu_2)((A_1-A_{2,s})w_0,
\overline{w_{1,s}})
e^{2i\tau\psi}d\sigma +\int_{\widetilde
\Gamma}(\nu_1+i\nu_2)((B_1-B_{2,s}) \widetilde w_0, \overline{\widetilde
w_{1,s}})e^{ -2i\tau\psi}d\sigma
\end{equation}
$$
+ o(\frac{1}{\root\of{\tau}}).
$$
Let $\Phi$ be given in Proposition \ref{Proposition -2}.
Then by (\ref{xoxo1u}),(\ref{iiii}) and the stationary phase argument,
the asymptotic formula holds:
\begin{equation}
\frak G(u_1,v)=\frac{1}{\root\of \tau}\sum_{x\in\mathcal
G}\{(\nu_1-i\nu_2)((A_1-A_{2,s})w_0,\overline{w_{1,s}})
e^{2i\tau\psi} +(\nu_1+i\nu_2)((B_1-B_{2,s}) \widetilde w_0,
\overline{\widetilde w_{1,s}})e^{-2i\tau\psi}\}(x)
\end{equation}
$$
+ o(\frac{1}{\root\of{\tau}})\quad\mbox{as}\,\,\tau\rightarrow
+\infty.
$$

Since for any $\widehat x$ one can find $\Phi$ such that $\widehat
x\in\mathcal G$ and $Im\Phi(\widehat x)\ne Im\Phi(x)$ for any $x \in
\mathcal G\setminus\{\widehat x\}$, we have
$$
((A_1-A_{2,s})w_0,\overline{w_{1,s}})=((B_1-B_{2,s}) \widetilde w_0,
\overline{\widetilde w_{1,s}})=0\quad \mbox{on}\,\,\Gamma_0.
$$
These equalities and Proposition \ref{nikita} imply (\ref{op!}).

Next we claim that
\begin{equation}\label{-15}
\frak G(e^{\tau\varphi}u_{-1}, v)=\frak G(u_{1},e^{-\tau\varphi} v_{-1})
=o(\frac 1\tau )\quad\mbox{as}\,\,\tau\rightarrow +\infty.
\end{equation}
Obviously, by (\ref{mimino11}) and Proposition \ref{vasya}, we see that
\begin{equation}
\frak G(e^{\tau\varphi}u_{-1},v- \mathcal V)=o(\frac 1\tau
)\quad\mbox{as}\,\,\tau\rightarrow +\infty.
\end{equation}

Let $\chi\in C_0^\infty(\Omega)$ satisfy
$\chi\vert_{\Omega\setminus \mathcal O_{\frac \epsilon 2}}=1.$
By (\ref{mimino11}), we have
\begin{eqnarray}\label{llg1}
\frak G(e^{\tau\varphi}u_{-1}, \mathcal V)=\frak
G(e^{\tau\varphi}u_{-1}, \chi\mathcal V)+o(\frac 1\tau) \nonumber\\
= \int_\Omega (2(A_1-A_{2,s})\partial_z(e^{\tau\varphi}u_{-1})+2(B_1-B_{2,s})
\partial_{\overline z}(e^{\tau\varphi}u_{-1}),
\chi\overline {\mathcal V})dx + o(\frac 1\tau) \nonumber\\
= \int_\Omega
(2(A_1-A_{2,s})\partial_z(e^{\tau\varphi}u_{-1}),
\overline{\chi\widetilde
w_{1,s}e^{\tau\overline\Phi}})+(2(B_1-B_{2,s})\partial_{\overline
z}(e^{\tau\varphi}u_{-1}),\overline{\chi
w_{1,s}e^{\tau\Phi}})dx+o(\frac 1\tau).
\end{eqnarray}
Let functions $w_4,w_5$ solve the equations $(-\partial_{\overline
z}+B^*_1)w_4=2(A_1-A_{2,s})^*{\widetilde w_{1,s}}$ and $(-\partial_{
z}+A^*_1)w_5=2(B_1-B_{2,s})^*{w_{1,s}}.$

Taking the scalar product of equation (\ref{zad2}) and the function
$w_5e^{\tau\Phi}+w_4e^{\tau\overline\Phi}$, after integration by parts we
obtain
\begin{equation}\label{llg}
\int_\Omega(2\partial_z
(e^{\tau\varphi}u_{-1})+A_1(e^{\tau\varphi}u_{-1}),
2(A_1-A_{2,s})^*\overline{\widetilde
w_{1,s}e^{\tau\overline\Phi}})
\end{equation}
$$
+(2\partial_{\overline z}
(e^{\tau\varphi}u_{-1})+B_1(e^{\tau\varphi}u_{-1}),2(B_1-B_{2,s})^*
\overline{w_{1,s}e^{\tau\Phi}})dx=o(\frac{1}{\tau}).
$$
By (\ref{llg1}) and (\ref{llg}), we obtain the first equality in
(\ref{-15}). The proof of the second equality in (\ref{-15}) is the
same.

By (\ref{zzz}), (\ref{op!}), (\ref{BIGzopa}), (\ref{-10}),
(\ref{-11}), (\ref{MAHA1}), (\ref{MAHA}) and (\ref{-15}), we have the
asymptotic formula:
\begin{eqnarray}\label{poso1}
 e^{2i\tau\psi(\widetilde
x)}(((B_1-B_{2,s})\partial_{\overline z}\overline\Phi b_+,
w_{1,s})_{L^2(\Omega)}+((B_1-B_{2,s})\partial_{\overline z}\overline\Phi
\widetilde w_0,\widetilde a_-)_{L^2(\Omega)})\\+ e^{-2i\tau\psi(\widetilde
x)}(((B_1-B_{2,s})\partial_{\overline z}\overline\Phi b_-,
w_{1,s})_{L^2(\Omega)}+((B_1-B_{2,s})\partial_{\overline z}\overline\Phi
\widetilde w_0,\widetilde
a_+)_{L^2(\Omega)})\nonumber\\
+ e^{2i\tau\psi(\widetilde x)}(((A_1-A_{2,s})\partial_z\Phi a_+,\widetilde
w_{1,s})_{L^2(\Omega)}+((A_1-A_{2,s})\partial_z\Phi w_0,\widetilde
b_-)_{L^2(\Omega)})\nonumber\\+ e^{-2i\tau\psi(\widetilde
x)}(((A_1-A_{2,s})\partial_z\Phi a_-, \widetilde
w_{1,s})_{L^2(\Omega)}+((A_1-A_{2,s})\partial_z\Phi w_0, \widetilde
b_+)_{L^2(\Omega)})\nonumber\\
-\pi\frac{ ({\mathcal Q}_+w_0,\overline w_1)e^{2i\tau\psi(\widetilde
x)+s}}{\tau\vert \det \psi''(\widetilde x)\vert^\frac 12} -\pi\frac{
({\mathcal Q}_-\widetilde w_0,\overline{\widetilde
w_1})e^{-2i\tau\psi(\widetilde x)+s}}{\tau\vert \det \psi''(\widetilde
x)\vert^\frac 12}+\mathcal P(\tau)+o(\frac 1\tau), \nonumber
\end{eqnarray}
where ${\mathcal Q}_+=2\partial_z(A_1-A_2)+B_2(A_1-A_{2})+(B_1-B_2)A_1
-(Q_1-Q_2)$ and
${\mathcal Q}_-=2\partial_{\overline
z}(B_1-B_2)+A_1(B_1-B_2)+(A_1-A_2)B_1-(Q_1-Q_2)$ and
\begin{eqnarray}
 \mathcal P(\tau)=-2\tau\int_\Omega (\mbox{\bf T}^*_{-
A_{2}^*}((A_1-A_{2,s})\partial_z\Phi   e^{s\eta} w_0) ,
\overline{q_{3}+\widetilde q_{3}/\tau})e^{2i\tau\psi}dx\nonumber\\
-\tau\int_\Omega(2\mbox{\bf
P}^*_{-B_{2}^*}((B_1-B_{2,s})\partial_{\overline z}\overline \Phi
e^{s\eta}\widetilde w_0), \overline
{q_{4}+\widetilde q_{4}/\tau })e^{-2i\tau\psi}dx\nonumber\\
-2\tau\int_\Omega e^{-2i\tau \psi}(q_2+\widetilde q_2/\tau ,\overline{
\mbox{\bf P}^*_{A_1}((A_1-A_{2,s})^*(\partial_{\overline
z}\overline\Phi{\widetilde w}_{1,s}))})dx\nonumber\\
-2\tau\int_\Omega (q_1+\widetilde q_1/\tau,\overline{\mbox{\bf
T}^*_{B_1}((B_1-B_{2,s})^*( {\partial_z\Phi} w_{1,s}))})e^{2i\tau
\psi}dx\nonumber
\end{eqnarray}

Observe that
$$
\mbox{\bf T}^*_{-A_{2}^*}((A_1-A_{2,s})e^{s\eta}w_0)+e^{s\eta}w_0\in Ker\,
\mbox{\bf T}^*_{-A_{2}^*}, \,\, \mbox{\bf
P}^*_{-B_{2}^*}((B_1-B_{2,s})e^{s\eta}\widetilde w_0)+e^{s\eta}
\widetilde w_0\in Ker\,
\mbox{\bf P}^*_{-B_{2}^*}
$$
and
$$
\mbox{\bf P}^*_{A_1}((A_1-A_{2,s})^*\widetilde w_{1,s})+\widetilde
w_{1,s}\in Ker\,\mbox{\bf P}^*_{A_1},\,\,\mbox{\bf
T}^*_{B_1}((B_1-B_{2,s})^* w_{1,s})+ w_{1,s}\in Ker\, \mbox{\bf
T}^*_{B_1}.
$$

Thanks to Proposition \ref{gandon} and above relations, there exist
functions $r_{1,s}\in Ker \,\mbox{\bf T}^*_{-A_{2}^*},r_{2,s}\in
Ker\, \mbox{\bf P}^*_{-B_{2}^*},r_{3,s}\in Ker\,\mbox{\bf
P}^*_{A_1},r_{4,s}\in Ker\, \mbox{\bf T}^*_{B_1}$ such that
\begin{eqnarray}
\mathcal P(\tau)=2\tau\int_\Omega (\partial_z\Phi  w_0,
\overline{e^{\tau(\overline\Phi-\Phi)}q_{3,s})}dx+\tau\int_\Omega
(r_{1,s},
\overline{e^{\tau(\overline\Phi-\Phi)}q_{3})}dx\nonumber\\
+2\tau\int_\Omega(\partial_{\overline z}\overline \Phi
\widetilde w_0, \overline
{e^{\tau(\Phi-\overline\Phi)}q_{4,s}})dx+\tau\int_\Omega
( r_{2,s}, \overline {e^{\tau(\Phi-\overline\Phi)}q_{4}})dx\nonumber\\
+2\tau\int_\Omega e^{\tau (\overline
\Phi-\Phi)}\partial_z\Phi(q_2+\widetilde q_2/\tau, \overline {\widetilde
w_{1,s}})dx+\tau\int_\Omega e^{\tau (\overline \Phi-\Phi)}(q_2,
\overline r_{3,s})dx\nonumber\\
+2\tau\int_\Omega e^{\tau ( \Phi-\overline
\Phi)}q_1,\overline{\partial_z\Phi} \overline{
w_{1,s}})dx+\tau\int_\Omega e^{\tau ( \Phi-\overline \Phi)}(q_1,\overline
r_{4,s})dx+o(\frac 1\tau)\quad\mbox\,\,{as}\,\,\tau\rightarrow
+\infty .
\end{eqnarray}
Integrating by parts in the above equality, we have
\begin{eqnarray}
\mathcal P(\tau)=-2\int_\Omega (w_0,
\overline{e^{\tau(\overline\Phi-\Phi)}\partial_z
q_{3,s})}dx+\tau\int_\Omega (r_{1,s},
\overline{e^{\tau(\overline\Phi-\Phi)}q_{3})}dx\nonumber\\
-2\int_\Omega( \widetilde w_0, \overline
{e^{\tau(\Phi-\overline\Phi)}\partial_{\overline z}q_{4,s}})dx
+ \tau\int_\Omega( r_{2,s}, \overline {e^{\tau(\Phi-\overline\Phi)}
q_{4}})dx\nonumber\\
+2\int_\Omega e^{\tau (\overline \Phi-\Phi)}(\partial_z q_2,
\overline {\widetilde w}_{1,s})dx+\tau\int_\Omega
e^{\tau (\overline \Phi-\Phi)}(q_2,\overline
r_{3,s})dx\nonumber\\
+2\int_\Omega e^{\tau ( \Phi-\overline \Phi)}(\partial_{\overline z}q_1,
\overline
w_{1,s})dx+\tau\int_\Omega e^{\tau ( \Phi-\overline \Phi)}(q_1,\overline
r_{4,s})dx+o(\frac 1\tau)\quad\mbox{as}\,\,\tau\rightarrow +\infty .
\end{eqnarray}
Applying the stationary phase argument, we obtain that there exists a
constant $\mathcal C_3$, independent of $\tau$, such that
\begin{eqnarray}\label{poso}
\mathcal P(\tau)= \frac{\mathcal C_3}{\tau}-2\pi\frac{ ({\mathcal
Q}_+w_0,\overline w_1)e^{2i\tau\psi(\widetilde x)+s}}{\tau\vert \det
\psi''(\widetilde x)\vert^\frac 12} -2\pi\frac{ ({\mathcal
Q}_-\widetilde w_0,\overline{\widetilde
w_1})e^{-2i\tau\psi(\widetilde x)+s}}{\tau\vert \det
\psi''(\widetilde x)\vert^\frac 12}\nonumber\\
+\frac{e^{2i\tau\psi(\widetilde x)}}{\tau \vert \det \psi''(\widetilde
x)\vert^\frac 12}(\frak D( \ell_4)+ \frak D(
\ell_2))+\frac{e^{-2i\tau\psi(\widetilde x)}}{\tau \vert \det
\psi''(\widetilde x)\vert^\frac 12}(\frak D( \ell_1)+\frak D( \ell_3)),
\end{eqnarray}
where $\ell_1=(q_1,\overline{ r_{4,s}}),\ell_2=(q_2,\overline{
r_{3,s}}), \ell_3=( r_{2,s}, \overline {q_{4}}), \ell_4=(r_{1,s},
\overline{q_{3}})$ and  for any smooth function $\ell(x)$ we set
$$\frak D( \ell)=\left(\partial_z (\frac{\ell_z(\widetilde x)(z-\widetilde
z)}{\partial_z\Phi})-\partial_{\overline z} (\frac{\ell_{\overline z}
(\widetilde x)(\overline z-\overline{\widetilde z})}
{\partial_{\overline z}\overline\Phi})+\frac
12\partial_z (\frac{\ell_{zz}(\widetilde x)(z-\widetilde
z)^2}{\partial_z\Phi})-\frac 12\partial_{\overline z} (\frac{\ell_{\overline
z\overline z}(\widetilde x)(\overline z-\overline{\widetilde z})^2}
{\partial_{\overline z}\overline\Phi})\right )(\widetilde x).
$$

Since $\psi(\widetilde x)\ne 0$, we obtain from (\ref{poso}) and
(\ref{poso1}):
\begin{eqnarray}\label{012}
 ((B_1-B_{2,s})\partial_{\overline z}\overline\Phi b_+,
w_{1,s})_{L^2(\Omega)}+((B_1-B_{2,s})\partial_{\overline z}\overline\Phi
\widetilde
w_0,\widetilde a_{-,s})_{L^2(\Omega)}\\
+ ((A_1-A_{2,s})\partial_z\Phi a_+,\widetilde
w_{1,s})_{L^2(\Omega)}+((A_1-A_2)\partial_z\Phi w_0,\widetilde
b_{-,s})_{L^2(\Omega)}\nonumber\\
-\pi\frac{ ({\mathcal Q}_+w_0,\overline w_1)e^s}{\vert \det \psi''(\widetilde
x)\vert^\frac 12} +\frac{\frak D( \ell_4)+ \frak D( \ell_2)}{\vert
\det \psi''(\widetilde x)\vert^\frac 12}=0\nonumber
\end{eqnarray}
and
\begin{eqnarray}\label{013}
((B_1-B_{2,s})\partial_{\overline z}\overline\Phi b_-,
w_{1,s})_{L^2(\Omega)}+((B_1-B_{2,s})\partial_{\overline z}\overline\Phi
\widetilde w_0,\widetilde
a_{+,s})_{L^2(\Omega)}\nonumber\\
+ ((A_1-A_{2,s})\partial_z\Phi a_-, \widetilde
w_{1,s})_{L^2(\Omega)}+((A_1-A_{2,s})\partial_z\Phi w_0, \widetilde
b_{+,s})_{L^2(\Omega)}\nonumber\\
-\pi\frac{ ({\mathcal Q}_-\widetilde w_0,\overline{\widetilde
w_1})e^s}{\vert \det \psi''(\widetilde x)\vert^\frac 12}+ \frac{\frak D(
\ell_1)+\frak D( \ell_3)}{\vert \det \psi''(\widetilde x)\vert^\frac
12}=0.\end{eqnarray}

Integrating by parts in (\ref{012}) and (\ref{013}), we obtain
\begin{eqnarray}\label{012A}
 ((\nu_1-i\nu_2)\partial_{\overline z}\overline\Phi b_+,
w_{1})_{L^2(\partial\Omega)}+((\nu_1-i\nu_2)\partial_{\overline
z}\overline\Phi \widetilde
w_0,\widetilde a_{-})_{L^2(\partial\Omega)}\\
+ ((\nu_1+i\nu_2)\partial_z\Phi a_+,\widetilde
w_{1})_{L^2(\partial\Omega)}+((\nu_1+i\nu_2)\partial_z\Phi
w_0,\widetilde
b_{-})_{L^2(\partial\Omega)}\nonumber\\
-\pi\frac{ ({\mathcal Q}_+w_0,\overline w_1)e^s}
{\vert \det \psi''(\widetilde
x)\vert^\frac 12} +\frac{\frak D( \ell_4)+ \frak D( \ell_2)}{\vert
\det \psi''(\widetilde x)\vert^\frac 12}=0               \nonumber
\end{eqnarray}
and
\begin{eqnarray}\label{013A}
((\nu_1-i\nu_2)\partial_{\overline z}\overline\Phi b_-,
w_{1})_{L^2(\partial\Omega)}+((\nu_1-i\nu_2)\partial_{\overline
z}\overline\Phi \widetilde w_0,\widetilde
a_{+})_{L^2(\partial\Omega)}\nonumber\\
+ ((\nu_1+i\nu_2)\partial_z\Phi a_-, \widetilde
w_{1})_{L^2(\partial\Omega)}+((\nu_1+i\nu_2)\partial_z\Phi w_0,
\widetilde
b_{+})_{L^2(\partial\Omega)}\nonumber\\
-\pi\frac{ ({\mathcal Q}_-\widetilde w_0,\overline{\widetilde
w_1})e^s}{\vert \det \psi''(\widetilde x)\vert^\frac 12}+ \frac{\frak D(
\ell_1)+\frak D( \ell_3)}{\vert \det \psi''(\widetilde x)\vert^\frac
12}=0.
\end{eqnarray}
Observe that
\begin{equation}\label{POpl}
\sum_{k=1}^4\vert \frak D( \ell_k )\vert\le C_4
\end{equation}
with the constant $C_4$ independent of $s.$  We prove this
inequality for $ \frak D( \ell_4).$ The proof  for remaining terms
is similar. By (\ref{BIGzopa}) the functions  $
(A_1-A_{2,s})e^{s\eta} w_0$ are bounded uniformly in the space
$H^1(\Omega)'.$ Hence, by (\ref{sonka1},the functions $\mbox{\bf
T}^*_{-A_2}(\partial_z\Phi(A_1-A_{2,s})e^{s\eta}
 w_0)$ are uniformly bounded in $L^2(\Omega)$. Then the functions
 $r_{1,s}$ are uniformly  bounded in $L^2(\Omega)$ and $Ker\,\mbox{\bf
 T}^*_{-A_2}.$ Therefore the functions $r_{1,s}$ are uniformly  bounded
 in $C^5(K)$ for any compact $K\subset\subset \Omega.$ Since $\ell_1=(r_{1,s},
\overline{q_{3}})$, the proof of (\ref{POpl}) is completed.

Passing to the limit in (\ref{012A}) and (\ref{013A}) as $s$ goes to
infinity, we obtain $({\mathcal Q}_+w_0,\overline w_1)(\widetilde x)=
({\mathcal Q}_-\widetilde w_0,\overline{\widetilde w_1})(\widetilde x)=0. $
These equalities and (\ref{xoxo1}) imply the equalities (\ref{A1}) and
(\ref{A2}) at point $\widetilde x$. According to Proposition
\ref{Proposition -1}, a point $\widetilde x$ can be chosen arbitrarily close
to any point of domain $\Omega$ after an
appropriate choice of the function $\Phi.$ The proof of the theorem
is completed.

$\blacksquare$

{\bf Acknowledgements.}
Most part of the paper has been written during the stay of
the first named author at Graduate School of Mathematical
Sciences of The University of Tokyo and he thanks the Global
COE Program ``The Research and Training Center for New Development
in Mathematics" for support of the visit to The University of Tokyo.

\end{document}